\documentclass[review]{elsarticle}

\usepackage{hyperref}

\usepackage{amsfonts,amsmath} 
\usepackage{algpseudocode}
\usepackage{algorithm}
\usepackage{multirow}
\usepackage{slashbox}
\usepackage{subfigure}
\usepackage{graphicx}

\journal{Computers \& Mathematics with Applications}

\newcommand{\R}{\mathbb{R}}
\newcommand{\N}{\mathbb{N}}
\newcommand{\mb}[1]{\mathbf{#1}}
\newcommand\mbx{\mb{x}}
\newcommand\mbk{\mb{k}}
\newcommand\mbi{\mb{i}}

\newcommand\mbF{\mb{F}}
\newcommand\mbG{\mb{G}}
\renewcommand\div[1]{{\rm div\,} {#1}}

\newcommand\e[1]{{\rm e}^{#1}}
\newcommand{\E}{\mathbb{E}}
\newcommand \Esp[1]{\E^Q\left[#1\right]}

\newcommand \dx[1]{{\rm d}#1}

\newcommand\Real[1]{\operatorname{Re}\left(#1\right)}

\begin{document}

\begin{frontmatter}

\title{A second order finite volume IMEX Runge-Kutta scheme for two dimensional PDEs in finance\tnoteref{mytitlenote}}
\tnotetext[mytitlenote]{M.J Castro research has been partially supported by the Spanish Government and FEDER through the coordinated Research project RTI2018-096064-B-C1, the Junta de Andaluc\'ia research project P18-RT-3163, the Junta de Andalucia-FEDER-University of M\'alaga research project UMA18-FEDERJA-16 and the University of M\'alaga. 
The other authors research has been partially supported by the Spanish MINECO under research project number PDI2019-108584RB-I00 and by the grant ED431G 2019/01 of CITIC, funded by Conseller\'ia de Educaci\'on, Universidade e Formaci\'on Profesional of Xunta de Galicia and FEDER.}

\author[auth1]{J. G. L\'opez-Salas}\ead{jose.lsalas@udc.es}
\address[auth1]{Department of Mathematics, Faculty of Informatics and CITIC, Campus Elvi\~na s/n, 15071-A Coru\~na (Spain)}

\author[auth1]{M. Su\'arez-Taboada}\ead{maria.suarez3@udc.es}

\author[auth2]{M. J. Castro\corref{cor1}}\ead{mjcastro@uma.es}
\cortext[cor1]{Corresponding author.}
\address[auth2]{Department of An\'alisis Matem\'atico, Facultad de Ciencias, University of M\'alaga, Campus de Teatinos s/n, M\'alaga, 29080-Andaluc\'ia (Spain)}

\author[auth1]{A. M. Ferreiro-Ferreiro}\ead{ana.fferreiro@udc.es}
\author[auth1]{J. A. Garc\'ia-Rodr\'iguez}\ead{jose.garcia.rodriguez@udc.es}

\begin{abstract}

In this article we present a novel and general methodology for building second order finite volume implicit-explicit (IMEX) numerical schemes for solving two dimensional financial parabolic PDEs with mixed derivatives. In particular, applications to basket and Heston models are presented. The obtained numerical schemes have excellent properties and are able to overcome the well-documented difficulties related with numerical approximations in the financial literature. The methods achieve true second order convergence with non-regular initial conditions. Besides, the IMEX time integrator allows to overcome the tiny time-step induced by the diffusive term in the explicit schemes, also providing very accurate and non-oscillatory approximations of the Greeks. Finally, in order to assess all the aforementioned good properties of the developed numerical schemes, we compute extremely accurate semi-analytic solutions using multi-dimensional Fourier cosine expansions. A novel technique to truncate the Fourier series for basket options is presented and it is efficiently implemented using multi-GPUs.
\end{abstract}


\end{frontmatter}


\section{Introduction}\label{sec1}

The goal of this article is to develop a novel methodology for pricing financial options. A general  framework for building second order  numerical schemes for financial advection-diffusion-reaction PDEs with mixed   derivatives is designed. 
The technique is based in the use of finite volume implicit-explicit (IMEX) Runge-Kutta numerical methods.
This approach is applied for developing a second order numerical method for solving basket and Heston models in two spatial dimensions. To our knowledge, this is the first time that finite volume IMEX Runge-Kutta numerical schemes have been successfully applied to the numerical solution of parabolic financial PDEs.
The interplay between finite volume methods and IMEX time integrators is an extremely powerful combination in mathematical finance.
On the one hand, finite volume methods are well suited to treat convective terms. In fact, it allows to build schemes with true second order convergence, even in the presence of non-regular initial or boundary conditions, which is a well-known difficulty in the financial literature. 
On the other hand,  IMEX time integrator makes possible to overcome the tiny time steps induced by the spatial semi-discretization of the diffusive terms.  This is also of paramount importance, otherwise the problem is intractable from the computational point of view.  On top of that, the presented numerical schemes provide excellent finite difference approximations of the Greeks. Actually, these approximations are free of the spurious oscillations which typically appear when no strongly stable time integrators, such as Crank–Nicolson, are considered. 
Last but not least, in order to assess the robustness, accuracy and order of convergence of the proposed finite volume IMEX numerical schemes, the multidimensional Fourier cosine method is implemented. A novel way to truncate the series expansion for basket options is presented. Additionally, an efficient parallel implementation in a multi-GPU environment allows to recover extraordinarily accurate option prices. 

Option pricing models are of great importance in the financial industry. Several types of stochastic models can be found in the literature for modelling the price of an option. These  models are formulated as stochastic differential equations (SDEs) or systems of SDEs. They are based on Brownian motions (or Wiener processes) that are used to simulate the evolution of the underlying assets and sometimes also their volatilities.

Option pricing methods can be mainly classified into two categories, deterministic and stochastic algorithms. The majority of the first ones deal with the solution of PDEs, while the second ones employ the Monte Carlo method based on the formulation of the model in terms of expectations. As we have said before, in this article we will focus on numerical methods to solve PDEs. In next paragraphs we point out the main drawbacks of Monte Carlo with respect to PDE based numerical methods.

Firstly, we note that Monte Carlo rate of convergence is very slow, of order $O(\frac{1}{\sqrt{S}})$ for all dimensions, $S$ being the number of simulations (see \cite{glasserman-03}). 
Secondly, if the risk neutral distribution of the underlying at the expiry of the option is known explicitly, there is no need to discretize the SDE to perform Monte Carlo simulation. In that situation, Monte Carlo will be faster than PDE methods, specially in high dimensions, since it does not suffer the curse of dimensionality. However, in several option pricing models this is not the case, so that the SDEs have to be discretized. In this procedure we must choose a suitable time step $\Delta t$ for the discretization of the SDE, so that the discretization error is of the same order of magnitude as the standard error of the simulation. Therefore, when increasing the number of simulations $S$, one must increase the number of time steps so that the time discretization error remains of the same order as the standard error of the simulation. Monte Carlo theory suggests that $\Delta t$ should be proportional to $\frac{1}{\sqrt{S}}$. Thus, in a simulation with $S$ paths we take $\Delta t = \kappa\frac{1}{\sqrt{S}} $ for some constant of proportionality $\kappa$. Summing up, in order to increase Monte Carlo accuracy one should need to increase the number of simulations which in turn requires a decrease in the discretization time step, and leads to a very slow convergence.

Additionally, in the context of options pricing, standard Monte Carlo methods have difficulties in providing accurate sensitivities (Greeks) of the option prices with respect to the underlying stocks or parameters, which are important for trading purposes. These Greeks are directly given by the discretization of the partial derivatives involved in the PDE or can be easily obtained by numerical differentiation formulas.
Also, early exercise features like in American or Bermudan options do not fit into the forward simulation paradigm of pure Monte Carlo and the combination with regression techniques turns out very costly (see \cite{gobetSalas} for details), whereas numerical methods for complementarity PDE problems are much more straightforward.

In this work we will focus in pricing options under financial models with two spatial variables. More precisely, the models we will be dealing with can be written as a 2D advection-diffusion-reaction PDE of the general form
\begin{equation}
\label{eq:AdvDiffPDE}
\dfrac{\partial u}{\partial t}
+a\dfrac{\partial u}{\partial x}
+b\dfrac{\partial u}{\partial y}
+c\dfrac{\partial^2 u}{\partial x^2}+
d\dfrac{\partial^2 u}{\partial y\partial x}+ 
e\dfrac{\partial^2 u}{\partial y^2}
+f=0,
\end{equation}
where $a,b,c,d,e$ are functions of $t,x$ and $y$, and $f = f(t,x,y,u)$. 

Standard numerical methods for solving parabolic PDEs in finance, of the type \eqref{eq:AdvDiffPDE}, are finite difference (FD) based, see \cite{Pironneau05,Duffy06,Wilmott06,tavellaRandall2000}, or finite elements (FE) based, see \cite{ACHDOU-PIRONNEAU-02,Pironneau05,VazquezNogueiras2006,VazquezNogueiras2006-2}.
Another possibility is to use Exponentially Fitted (EF) finite differences schemes (see \cite{Duffy06}) or Alternate Directions (ADI) finite differences methods (see \cite{Foulon10,HENDRICKS2017175}).

Nevertheless, standard FD or FE numerical schemes for financial PDEs have to face also several mathematical challenges. The first one is due to the lack of smoothness of the terminal conditions imposed by most of the traded option contracts. A financial option is a contract that gives the holder the right to trade in the future at a previously agreed exercise price. However, the option holder is under no obligation and will only exercise if it is worth it. This optionality in the exercise is  represented mathematically in the payoff function of the contract, commonly by means of the $\max$ function, i.e $\max(\cdot,0)$. As a consequence, terminal conditions of PDE models, and their derivatives, are commonly non-smooth. This non-smooth data originates noise in the numerical solution. As a result, the expected rate of convergence is seriously degraded. Besides, the numerical computation of derivatives of the numerical solution originates poor estimates of the Greeks with oscillations. Several strategies were designed in the literature in order to overcome these issues. For instance, smoothing methods for the terminal conditions were proposed. Heston and Zhou \cite{hestonZhou2000} replace mesh nodal values with an average value computed over a close space. This technique is natural in finite volume methods. In \cite{tavellaRandall2000}, Tavella and Randall proposed to work with meshes such that the strike price is halfway between grid points. Another method consists of projecting the terminal conditions onto a set of basis functions, see \cite{Rannacher1984}. Nevertheless, in general a smoothing technique must be combined with thorough time discretizations in order to recover second order convergence for option pricing problems with discontinuous payoff functions (or its derivatives). Along this line, Rannacher \cite{Rannacher, Rannacher1984} proposed a numerical scheme in which a finite number of fully implicit time steps are performed before the Crank-Nicolson time discretization. This approach, known as Rannacher timestepping, along with a suitable smoothing technique of the terminal condition, preserves the quadratic convergence rate of the Crank-Nicolson scheme (see  \cite{gilesCarter,PooleyVetzalForsyth} for more details).

Another mathematical challenge in the numerical solution of  this kind of advection-diffusion PDEs arises when the advection term becomes larger than the diffusion one, or when we have degenerated diffusion. Under this situation, instabilities show up, because the problem becomes more hyperbolic. In order to face these problems, several techniques have been introduced. One way to overcome these instability phenomena is to avoid centered schemes and to consider upwind discretizations of advection terms, thus taking information upstream and not downstream. We refer the reader to \cite{Pironneau05} for further details in the mathematical finance field. Although upwind schemes are stable, they are less accurate than the centered ones.
Another approach to face instabilities in convection-dominated diffusion problems is the semi-Lagrangian time stepping method, that integrate the equations by searching backward in time the position of the flow. Interpolation is needed in this tracing procedure along characteristic curves. Several works in the literature combine the method of characteristics with FE or FD procedures in option pricing problems. In \cite{ForsythCharact06}, Forsyth et al. price American Asian options under jump diffusion. In \cite{VazquezNogueiras2006} the authors price Asian options under a semi-Lagrangian finite element method. In \cite{Ferreiro13}, Ferreiro et al. present a semi-Lagrangian finite difference method for pricing business companies. The main disadvantages of semi-Lagrangian methods is the difficulty to build high order numerical schemes. In fact, as stated in \cite{ForsythCharact06}, due to the non-smoothness of the payoff, ``\textit{even if second order time stepping methods are used, observed convergence as the mesh and time step is refined occurs at a sub-second order rate}''. Besides, the computational cost of this strategy can be high due to the backward search in time. In fact, at each time step an initial value problem for each point of the spatial mesh has to be solved, and then the position of the point has to be traced backward in time. Finally, interpolation with high enough order has to be performed. Another approach to properly solve convection-dominated diffusion problems is to consider finite volume discretization methods. Finite volume methods naturally handle non smooth payoffs and the cases where the PDE becomes diffusion degenerated. Seminal works on applying finite volume method to option pricing problems are due to  Forsyth and Zvan, see  \cite{ForsythZvan01, Zvan-Forsyth-97}, where vertex type finite volume methods were applied for problems in two spatial dimensions. More recently, in \cite{mathias13}, finite volume methods were proposed for convection dominated pricing problems. More precisely, the authors applied the Kurganov-Tadmor scheme proposed in \cite{NT1990,KT2000}  to option pricing problems written in conservative form.
Later, in \cite{Tadmor2ndOrder2019} the authors propose a second order improvement of \cite{mathias13} with appropriate time methods and slope limiters. Finally in \cite{Tadmor2019} the authors apply a third order Kurganov-Levy scheme presented in \cite{KL2000} based in the CWENO reconstructions presented in \cite{LPR1999}. All these works only deal with Black-Scholes PDEs in dimension one or Asian PDE problems. Therefore, they do not solve general two-dimensional pricing problems with mixed derivatives. In these works, the ordinary differential equations (ODEs) obtained after semi-discretizing in space are stiff. In fact, the authors use explicit-schemes in time. The main disadvantage of these approaches is that time steps have to be very small to maintain stability. Typically, Von Neumann stability analysis demands $\Delta t \leq \frac{(\Delta x)^2}{2 \eta}$, being $\eta$ the diffusion velocity.

Pareschi and Russo proposed in \cite{Russo05} implicit-explicit (IMEX) Runge-Kutta method for general hyperbolic systems of conservation laws with stiff relaxation terms. The idea is to apply an implicit discretization to the source terms and an explicit one to the nonstiff term. In \cite{toivanen,frutos}, the IMEX time marching scheme was used in mathematical finance only combined with finite differences and finite elements, respectively.

All in all, our goal is to first semi-discretize in space 2D Black-Scholes equations (having mixed derivatives) with second-order finite volume methods, thus properly treating convection terms and non-smooth payoffs. Later, we propose to integrate in time the resulting system of stiff ODEs by means of the second-order IMEX Runge-Kutta time marching scheme. In this way, we will apply an implicit discretization to the diffusion (stiff) part and an explicit one to the convection and source terms (non stiff). As a result, the time step restriction will depend only on the stability of the degenerate diffusion PDE, i.e, $\Delta t \leq \frac{\Delta x}{\alpha}$, being $\alpha$ the convection velocity. The here proposed numerical methods will provide also good approximations of the Greeks (without oscillations).

The organization of this paper is as follows. In Section \ref{sec:PDE_Models} we overview a couple of two-dimensional Black-Scholes PDE problems with huge interest in the financial industry.
Section \ref{sec:numerical-schemes} presents the proposed finite volume IMEX Runge-Kutta scheme for general two-dimensional convection diffusion PDEs written in conservative form. Section \ref{2dcos} is devoted to an alternative Fourier numerical method designed to solve the PDE problems presented in Section \ref{sec:PDE_Models}. In order to obtain a high-accurate reference solution to validate the here proposed finite volume methods, a novel and very efficient multi-GPU implementation of the Fourier method is carried out. In Section \ref{sec:numerical-results} numerical experiments are performed to price several options under the Black-Scholes and Heston models. A thorough inspection of the accuracy and order of convergence of the developed numerical methods is detailed. Finally, conclusions are drawn in Section \ref{conclusions}.

\paragraph{Notations of the paper} 
\begin{enumerate}[(i)]
 \item $\ln(x)$ stands for the natural logarithm of $x\in \mathbb{R}_+$.
 \item {$\Real{z}$ denotes the real part of the complex $z\in \mathbb{C}$; we set $\mb{i}:=\sqrt{-1}$.}
 \item $N$ is the standard normal cumulative distribution function.
 \item $\mb{x}^\mathsf{T}$ stands for the transpose of the vector $\mb{x}=(x_1,\dots, x_n)\in \R^n$.
 \end{enumerate}

 \section{Option pricing PDE models \label{sec:PDE_Models}}

The simplest option is a contract giving its owner the right (not the obligation) to buy (call option) or sell (put option) a particular asset at a fixed price at a specified time in the future. Therefore, the option contract has to specify three terms. The first one is the underlying asset; the price of this asset will be denoted by $s_t$. The second one is the fixed price, the so-called exercise price or strike price; it will be denoted by $K$. Finally, the given date is termed the maturity or the expiration date; it is often denoted by $T$. If the right to exercise the option is only permitted at maturity, the option is called European. This kind of options are basic derivatives, and are often called plain vanilla options. Option pricing is one of the many problems of financial mathematics. The price of the option at maturity is straightforward. If $s_T$ is the price of the underlying asset at maturity, the payoff of a call option is $\max(s_T-K,0)$, whilst $\max(K-s_T,0)$ is the benefit of the option's holder at the option's expiry. Nevertheless, determining the price that the buyer of the option has to pay to the seller at the initial time $t$, the so-called premium, is a challenging mathematical problem because $s_T$ is not known at time $t$. We will use $u$ to mean the value of the option. It is clear that $u$ is a function of the value of the underlying asset $s$ at time $t$. The goal of option pricing PDE models is to answer the following question: what is the market price $u(s,t)$ of the option at time $t$, $0\leq t < T$?

The theory of option pricing started in the seventies with the celebrated works of Black and Scholes \cite{BS73} and Merton \cite{Merton73}. The Black-Scholes model assumes that the price of the underlying asset is the solution of the following stochastic differential equation $$\frac{\dx{s_t}}{s_t} = r \dx{t} + \sigma \dx{W_t},$$ where $W$ is a standard Brownian motion, $r \in \R$ is the interest rate and $\sigma\in\R_+$ is the volatility of the underlying asset. Therefore, the vanilla call option price is given by the following expectation under the risk-neutral measure $$u(s_t,t) = \e{-r(T-t)}\mathbb{E}[\max(s_T-K,0)].$$ Besides $u$ solves the backward PDE $$\dfrac{\partial u}{\partial t} + \frac12 \sigma^2 s^2 \frac{\partial^2 u}{\partial s^2} + r S \frac{\partial u}{\partial s} - r u = 0,\quad t\in[0,T).$$ In this work we will work on deterministic methods for pricing options, thus proposing appropriate numerical methods for solving parabolic advection–diffusion PDEs. However, we will deal with more complicated contracts and models than the one described above, so that the resulting PDE becomes two-dimensional in space. More precisely, we are interested in pricing options on many underlyings, the so-called basket options, and also in pricing vanilla options under stochastic volatility models. 

On the one hand, an arithmetic basket call option has a payoff of  $$\max(c_1 s_{1T} + c_2 s_{2T} - K,0), \quad c_1,c_2\in\R.$$ The Black-Scholes model results in closed form solutions for pricing vanilla options, due to the fact that the underlying asset price at fixed time follows a lognormal distribution. Nevertheless, considering the Black-Scholes model for several underlying assets does not lead to closed form formulas for the price of a basket option. The obstacle is mainly the fact that the distribution of a weighted average of lognormals is not known. Thus, numerical methods are required to approximate  the price of arithmetic basket options.

On the other hand, the Black-Scholes model assumes that $\sigma$, the volatility of the underlying, is constant. However, this model cannot recover the market volatility smile and skew, which are behaviours of the market implied volatility surface (see \cite{gatheral} and references therein for details). In order to overcome this limitation of constant volatility models, stochastic volatility models were proposed. They change our model for the evolution of the underlying asset by adding a new variable for the volatility which satisfy some SDE. The most popular and renowned of all stochastic volatility models is the Heston model proposed in \cite{Heston93}. Pricing vanilla options under this model also requires numerical methods if we do not want to depend on quasi-closed form solutions.

In the next two subsections the problems we want to solve in this article are stated explicitly.

\subsection{Options on a basket of two assets\label{sec:Basket2D}}

We consider a basket of two assets, whose prices are $s_1$ and $s_2$. Under the Black-Scholes model, the prices of the underlying assets follow the following system of stochastic differential equations:
\begin{align} \label{eq:sdesBS2d}
  &\frac{\dx{s_{1t}}}{s_{1t}} =   (r-q_1) \dx{t} + \sigma_1 \dx{W_{1t}}, \quad s_{10} \mbox{ known}, \nonumber\\
  &\frac{\dx{s_{2t}}}{s_{2t}} =   (r-q_2) \dx{t} + \sigma_2\dx{W_{2t}}, \quad s_{20} \mbox{ known},
\end{align}
where $W$ is a two dimensional correlated Brownian motion. We set 
\begin{equation}\label{eq:brownianosCorrelados}W_{1t} = \bar{W}_{1t}, \quad  W_{2t} = \rho \bar{W}_{1t} + \sqrt{1-\rho^2} \bar{W}_{2t},\end{equation}
where $(\bar{W}_{1},\bar{W}_{2})$ is a two dimensional standard Brownian motion and $\rho \in (-1,1)$ is the constant correlation parameter. Besides, $\sigma_1, \sigma_2 \in \R_+$ are the market volatilities of the assets $s_1$ and $s_2$, respectively. Finally, $r$ is the interest rate, and $q_1$ and $q_2$ are the dividend yields of the assets.
The price of an option with maturity $T$ and payoff function $u_T(s_1,s_2)$ is $$u(s_{1t},s_{2t},t) = \e{-r(T-t)}\Esp{u_T(s_{1T},s_{2T})},$$ where $Q$ is the selected martingale measure.
By applying the two-dimensional Itô formula, we can find the PDE for the price of the option $u(s_1,s_2,t)$:
\begin{align}\label{eq:PdeBasket2DBackward}
    & \dfrac{\partial u}{\partial t}  + \frac12\sigma_1^2s_1^2\frac{\partial^2 u}{\partial s_1^2}+\frac12\sigma_2^2s_2^2\frac{\partial^2 u}{\partial s_2^2} + \rho\sigma_1\sigma_2s_1s_2\frac{\partial^2 u}{\partial s_1\partial s_2} \nonumber\\
    &\qquad +(r-q_1)s_1 \dfrac{\partial u}{\partial s_1}+(r-q_2)s_2\dfrac{\partial u}{\partial s_2}-ru=0,\qquad t\in[0,T), \nonumber\\
    & u(s_1,s_2,T) = u_T(s_1,s_2), \qquad s_1,s_2>0.
\end{align}
For convenience, a time reversal $\tau = T-t$ is performed. Doing so, the backward PDE \eqref{eq:PdeBasket2DBackward} is then transformed to the following forward PDE
\begin{align}\label{eq:PdeBasket2DForward}
    & \dfrac{\partial u}{\partial \tau }  - \frac12\sigma_1^2s_1^2\frac{\partial^2 u}{\partial s_1^2}-\frac12\sigma_2^2s_2^2\frac{\partial^2 u}{\partial s_2^2} - \rho\sigma_1\sigma_2s_1s_2\frac{\partial^2 u}{\partial s_1\partial s_2} \nonumber\\
    &\qquad -(r-q_1)s_1 \dfrac{\partial u}{\partial s_1}-(r-q_2)s_2\dfrac{\partial u}{\partial s_2}+ru=0,\qquad \tau\in(0,T], \nonumber \\
    & u(s_1,s_2,0) = u_0(s_1,s_2), \qquad s_1,s_2>0.
\end{align}
In the sequel, the $\tau$ notation is dropped for simplicity and the forward time is again written as $t$.
In this work we consider the following payoff function for the arithmetic basket call option
\begin{equation} \label{eq:arithmeticPayoff}
u_0(s_1,s_2)=\max\left(\dfrac{1}{2}(s_1+s_2)-K,0\right), 
\end{equation}
where $K$ is the fixed strike price.
As previously said, it is not possible to compute the value of an arithmetic basket call option analytically. We will compare the results of the proposed finite volume method with reference values calculated with another method, namely the COS method.

In order to compute a numerical approximation of $u$, we need to truncate the domain in the variables $s_1$ and $s_2$, so $u$ will be computed for $s_1,s_2\in(0,\bar{S})$, with $\bar{S}$ large enough. Besides, boundary conditions have to be imposed at the artificial boundaries $s_1=s_2=\bar{S}$. The PDE \eqref{eq:PdeBasket2DForward} is degenerate on the boundaries $s_1,s_2=0$ and is reduced into a Black-Scholes 1D PDE.
For the arithmetic basket call option with payoff function \eqref{eq:arithmeticPayoff} the following Dirichlet and Neumann boundary conditions can be used:
\begin{align*}
    u(0,s_2,t) & = C\left(\frac12 s_2, t, \sigma_1, q_1\right), & \frac{\partial u}{\partial s_1}(\bar{S},s_2,t) &= \frac12,\\
    u(s_1,0,t) & = C\left(\frac12 s_1,t, \sigma_2, q_2\right) , & \frac{\partial u}{\partial s_2}(s_1,\bar{S},t) &= \frac12,
\end{align*}
where $C$ is the Black-Scholes 1D price of the vanilla call European option, given by
\begin{align*}
  &C(S,t,\sigma,q)=S \e{-q t} N(d_1) - K \e{-r t} N(d_2),\\
  & d_1 = \frac{\ln(\frac{S}{K}) + (r-q + \frac{\sigma^2}{2})t }{\sigma \sqrt{t}}, \quad d_2 = d_1-\sigma\sqrt{t}.
\end{align*}

\subsection{Heston model \label{sec:HestonModel}}

In the Heston model \cite{Heston93}, there is an underlying asset $s$ whose volatility is a stochastic process driven by a second Brownian motion:
\begin{align}\label{SDE-System-Heston}
&\dx{s_t}=(r-q)s_t \dx{t} +\sqrt{v_t}s_t \dx{W_{1t}},\quad s_0 \mbox{ known,}\nonumber\\
&\dx{v_t}=\kappa(\theta-v_t)\dx{t} +\sigma\sqrt{v_t} \dx{W_{2t}}, \quad v_0 \mbox{ known},
\end{align}
where $r, q, \kappa, \theta, \sigma, s_0, v_0$ are constant parameters in $\R_+$ and $W$ is a two dimensional correlated Brownian motion, settled as in \eqref{eq:brownianosCorrelados}. Moreover, $r$ is the fixed interest rate, $q$ is the constant dividend yield, $\kappa$ is the mean-reversion speed for the variance, $\theta$ is the mean reversion level for the variance, $\sigma$ is the volatility of the variance (the so-called volatility of the volatility) and $v_0$ is  the initial level of the variance. The process $v$ represents the variance of $s$ and its stochastic differential equation is a version of the square root process described by Cox, Ingersoll and Ross \cite{coxIngersollRoss1985}.

The price of a derivative with payoff function $u_T(s)$ is the solution of the following backward PDE:
\begin{align}
\label{eq:HestonPDEBackward}
& \dfrac{\partial u}{\partial t}+
\dfrac{1}{2}s^2v\dfrac{\partial^2 u}{\partial s^2}+
\rho\sigma s v\dfrac{\partial^2 u}{\partial v\partial s}+
\dfrac{1}{2}\sigma^2 v \dfrac{\partial^2 u}{\partial v^2}\nonumber \\
& \qquad +(r-q)s \dfrac{\partial u}{\partial s}+
\kappa(\theta-v) \dfrac{\partial u}{\partial v}-
ru=0, \qquad t \in [0,T), \nonumber \\
&u(s,v,T) = u_T(s), \qquad s>0.
\end{align}
Once more, if we perform a time reversal $\tau = T-t$, the backward PDE \eqref{eq:HestonPDEBackward} is transformed to the following forward PDE
\begin{align}
\label{eq:HestonPDEForward}
& \dfrac{\partial u}{\partial \tau}-
\dfrac{1}{2}s^2v\dfrac{\partial^2 u}{\partial s^2}-
\rho\sigma s v\dfrac{\partial^2 u}{\partial v\partial s}-
\dfrac{1}{2}\sigma^2 v \dfrac{\partial^2 u}{\partial v^2} \nonumber \\
& \qquad -(r-q)s \dfrac{\partial u}{\partial s}-
\kappa(\theta-v) \dfrac{\partial u}{\partial v}+
ru=0, \qquad \tau \in (0,T], \nonumber \\
&u(s,v,0) = u_0(s), \qquad s>0.
\end{align}
As before, the forward time $\tau$ will be written as $t$.
The payoff function of an European call option is
\begin{equation} \label{eq:hestonPayoff}
u_0(s) = \max(s-K,0),
\end{equation}
where $K$ is the strike. Since the natural logarithmic characteristic function of $s$, $\mathbb{E}[\e{\mbi \xi \ln s_T}]$, is known explicitly, the price of the call option can be also computed by means of Fourier methods.

The finite volume numerical approximation of $u$ will be computed in a truncated domain, $s\in(0,\bar{S})$ and $v\in (0,\bar{V})$. In this work, the following boundary conditions were imposed $$u(0,v,t) = 0,\quad \frac{\partial^2 u}{\partial s^2}(\bar{S}, v, t) = 0, \quad \frac{\partial^2 u}{\partial v^2}(s, \bar{V}, t) = 0.$$ At the boundary $v=0$ no condition is specified. If the so-called Feller condition is satisfied, $2 \kappa \theta > \sigma^2$, it follows that this is and outflow boundary.

\section{Finite volume IMEX Runge-Kutta numerical method \label{sec:numerical-schemes}}

Equation \eqref{eq:AdvDiffPDE}  can be written in conservative form in the  following compact way:

\begin{equation} \label{SistCons}
    \dfrac{\partial}{\partial t}u(x,y,t)  + \dfrac{\partial f_1}{\partial x} (u)+\dfrac{\partial f_2}{\partial y}  (u)=\dfrac{\partial g_1}{\partial x} (u_x,u_y)+\dfrac{\partial g_2}{\partial y} (u_x,u_y)+h(u),
\end{equation}
or
\begin{equation} \label{SistConsCompacto}
    \dfrac{\partial}{\partial t}u(x,y,t)  + \div{\mbF}(u)=\div{\mbG}(\nabla u)+h(u),
\end{equation}
with $\mbF(u)=(f_1(u),f_2(u))$ and $\mbG(\nabla u)=(g_1(\nabla u),g_2(\nabla u))$. More precisely, $\mbF:\R^3\rightarrow\R^2$,  $\mbG:\R^4\rightarrow\R^2$ and $h:\R^3\rightarrow\R$ are functions of $x,y,u,\nabla u$, although in order to keep the notation simple some of these dependencies were omitted.

The numerical solution of equation \eqref{SistConsCompacto} using a finite volume semi-implicit scheme  is difficult because of the presence of the diffusive part. In this article we have considered the implicit-explicit (IMEX) Runge-Kutta time discretization numerical scheme proposed in \cite{Russo05}. This time-marching scheme play a major rule in the treatment of stiff systems of differential equations in the form
\begin{equation}\label{eq:stiffedo}
\dfrac{\partial U}{\partial t}  +  E (U)=I (U),
\end{equation}
where $U=U(t)\in \R^N$ and $E,I:\R^N\to \R^N$, being $E$ the non-stiff term and $I$  the stiff part. Both parts will be handled simultaneously with the same IMEX solver. 

The rest of the section is organized as follows. In Section \ref{sec:spatialDiscretization} we describe the space discretization obtained by finite volume schemes. Section \ref{sec:timeDiscretization} is devoted to the IMEX Runge-Kutta scheme applied to the stiff system of differential equations \eqref{eq:stiffedo} obtained by the finite volume space discretization.

\subsection{Space discretization} \label{sec:spatialDiscretization}

\begin{figure}
\centering
\includegraphics[width=8cm]{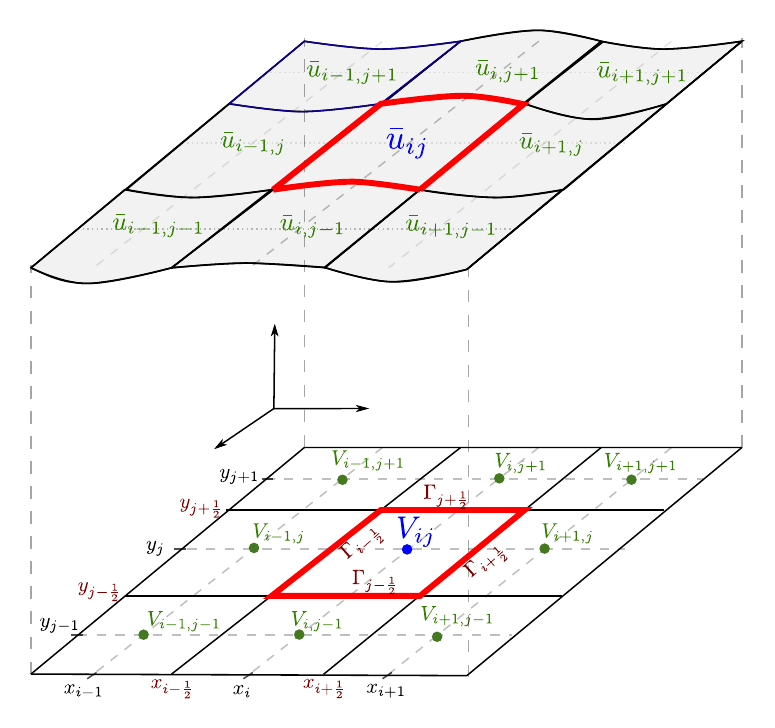}
\caption{Finite volume stencil.}\label{fig:Stencil}
\end{figure}

Let $\Delta x$, $\Delta y$ be the mesh length in the $X,Y$ directions respectively.
We define the grid points $x_i=i\Delta x$, $y_j=j\Delta y$, $i,j \in \mathbb{Z}$. Let $x_{i+1/2} = x_i + \frac12\Delta x$ and $y_{j+1/2} = y_j + \frac12\Delta y$. We consider rectangular finite volumes $V_{ij}= [x_{i -1/2} , x_{i +1/2} ]\times[y_{j -1/2} , y_{j +1/2} ] $, where $(x_i,y_j)$ is the center of the finite volume $V_{ij}$. The area of $V_{ij}$ will be denoted as $\lvert V_{ij}\rvert$, i.e.  $\lvert V_{ij}\rvert = \Delta x\Delta y$. Besides, $\Gamma_{ij}$ represents the boundary of $V_{ij}$. Let $\displaystyle\bar{u}_{ij}=\frac{1}{\lvert V_{ij}\rvert}\int_{V_{ij}} u(x,y,t) \dx{x}\dx{y}$ be the volume average. The here described finite volume grid is sketched in Figure \ref{fig:Stencil}.

Integrating equation \eqref{SistConsCompacto} in space on $V_{ij} $ and dividing by $\lvert V_{ij} \rvert$ we obtain the semi-discrete equation
\begin{align}
\label{eq:integrada_vols}
\dfrac{\dx{\bar{u}_{ij}}}{\dx{t}}
 =&-\dfrac{1}{ \lvert V_{ij}\rvert} \int_{V_{ij}}  \div{\mbF}(u) \,\dx{x}\dx{y} +\dfrac{1}{\lvert V_{ij} \rvert} \int_{V_{ij}}  \div{\mbG}(u_x,u_y) \,\dx{x}\dx{y}
\nonumber \\
&+\dfrac{1}{ \lvert V_{ij} \rvert} \int_{V_{ij}} h(u) \, \dx{x}\dx{y}.
\end{align}
Applying the divergence theorem we can rewrite the first two volume integrals of \eqref{eq:integrada_vols} as line integrals
\begin{align}
\label{eq:integragada_lineIntegrals}
\dfrac{\dx{ \bar{u}_{ij}}}{\dx{t}}
 =&-\dfrac{1}{ \lvert V_{ij} \rvert} \oint_{\Gamma_{ij}}  {\mbF}(u)\cdot{\mb{ n} } \,d\gamma 
+\dfrac{1}{ \lvert V_{ij} \rvert}  \oint_{\Gamma_{ij}}  {\mbG}(u_x,u_y)\cdot{\mb{ n} } \,d\gamma \nonumber \\
&+\dfrac{1}{ \lvert V_{ij} \rvert} \int_{V_{ij}} h(u) \, \dx{x}\dx{y}.
\end{align}
The unknowns of the problem are the volume averages $\{\bar{u}_{ij}(t)\}$. Therefore, in order to convert \eqref{eq:integragada_lineIntegrals} into a numerical scheme, we have to
approximate on the right hand side of this equation with functions of $\{\bar{u}_{ij}(t)\}$. In the following two content blocks, the numerical treatment of the source and advection terms (Section \ref{sec:numSemiDiscExplicit}) and the diffusion ones (Section \ref{sec:numSemiDiscImplicit}) is discussed. The source and convective terms will be treated explicitly by IMEX Runge-Kutta, while the diffusion part will be managed implicitly in time.

\subsubsection{Explicit part: advection-reaction space discretization} \label{sec:numSemiDiscExplicit}

For the advective part one readily gets
\begin{align}
 \oint_{\Gamma_{ij}} {\mbF}(u)\cdot{\mb{ n} } \,d\gamma=& \int_{\Gamma_{i+1/2},j} f_1(u) \dx{\gamma_j} + \int_{\Gamma_{i-1/2,j}} f_1(u) \dx{\gamma_j}  + \notag\\
 & \int_{\Gamma_{i,j+1/2}} f_2(u) ) \dx{\gamma_i} + \int_{\Gamma_{i,j-1/2}} f_2(u )  \dx{\gamma_i}, \label{eq:explicitPartAdvection}
\end{align} 
where for second order schemes each line integral can be approximated by means of  midpoint quadrature rule as follows
\begin{align*}
   \int_{\Gamma_{i\pm 1/2,j}} f_1(u) \dx{\gamma_j}&\approx \pm \Delta y f_1(u(x_{i\pm 1/2},y_j)),\\
   \int_{\Gamma_{i,j\pm1/2}} f_2(u) ) \dx{\gamma_i} &\approx  \pm \Delta x f_2(u(x_i,y_{j\pm1/2})).
\end{align*}

At this point the unknown function $u(x, y, t)$ is reconstructed by a piecewise polynomial using the volume averages $\{\bar{u}_{ij} (t)\}$. More precisely, starting from $\{\bar{u}_{ij}\}$, we compute a piecewise polynomial reconstruction
$$
\mathcal{R}(x,y) = \sum_{i,j}
 P_{ij} (x,y) \mathbf{1}_{ij} (x,y),
$$
where $P_{ij}$ is a suitable polynomial satisfying some accuracy and non oscillatory property, and $\mathbf{1}_{ij}$ is the indicator function of the volume $V_{ij}$. 
Second order numerical schemes can be obtained by means of piecewise linear polynomials, although higher order schemes can be obtained by polynomials of higher order. Here we consider the natural extension of MUSCL reconstruction to 2D Cartesian grids (see \cite{BramVanLeer}).

The flux functions at the midpoints of the boundaries of the volumes can be computed by using a
suitable numerical flux function, consistent with the analytical flux.
Hereafter we detail the approximation of $f_1 (u(x_{i +1/2},y_j))$, being the other fluxes approximated in the same way:
$$
f_1 (u(x_{i +1/2},y_j )) \approx \mathcal{F}_1 (u_{i+1/2,j}^- , u_{i +1/2,j}^+ ).
$$
The values $u^\pm_{i +1/2,j}$ are obtained from the reconstruction as
$$
u^\pm_{i +1/2,j}
 =
 \lim_{x\to x^\pm_{i +
 1/2}}
 \mathcal{R}(x,j\Delta y),
 $$
 with $x$ in a normal line to the boundary $\Gamma_{i+1/2,j}$. For example the left reconstructed value at the edge $\Gamma_{i+1/2,j}$ is
$$
u^-_{i + 1/2,j}= \bar{u}_{i,j}+\dfrac{\Delta x}{2} u'_{i,j},
$$
where the slope $u'_{i,j}$ is a first order approximation of the space derivative in the $X$ direction of $u(x,y,t)$ at point $(x_i, y_j)$ at every time $t$.
This slope must satisfy the TVD property and thus we must use slope limiters. In our case we use the minmod limiter, where the slope is given by
$$
u'_{ij}=\dfrac{1}{\Delta x} {\rm minmod} ( \bar{u}_{i,j}-\bar{u}_{i-1,j},\bar{u}_{i+1,j}-\bar{u}_{i,j}),
$$
$$
{\rm minmod}(a,b)=
\begin{cases}
\min (a,b) & \text{if } a,b>0,\\
\max(a,b) & \text{if } a,b<0,\\
0 & \text{otherwise\,. } 
\end{cases}
$$
All the calculations performed in this article have been carried out using the CIR numerical flux
$$
\mathcal{F}_1(u^-_{i +1/2,j},u^+_{i +1/2,j})=\dfrac{1}{2}(f_1(u^-_{i +1/2,j})+f_1(u^+_{i +1/2,j}))- \dfrac{1}{2}\alpha(u^+_{i +1/2,j}-u^-_{i +1/2,j}),
$$
where $\alpha=\left\lvert f'_1\left(\dfrac{u^-_{i +1/2,j}+u^+_{i +1/2,j}}{2}\right)\right\rvert$.

Therefore, the line integral \eqref{eq:explicitPartAdvection} of the convective term, is finally approximated as
\begin{align}
\oint_{\Gamma_{ij}}  {\mbF}(u)\cdot{\mb{ n} } \,d\gamma \approx& \Delta y \left( \mathcal{F}_1 (u^-_{i +1/2,j} , u^+_{i+1/2,j} ) -\mathcal{F}_1 (u^-_{ i -1/2,j} , u^+_{i-1/2,j} ) \right)+\notag \\  
 & \Delta x \left ( \mathcal{F}_2 (u^-_{i, j +1/2},u^+_{i,j+1/2} ) -\mathcal{F}_2 (u^-_{i,j -1/2},u^+_{i,j-1/2} ) \right),
\label{eq:semidiscretospaceConvection} 
\end{align}
being $\mathcal{F}_1$ and $\mathcal{F}_2$ the numerical fluxes of the physical flux functions $f_1$ and $f_2$, respectively.

Finally, the volume integral of the source term can be discretized using the midpoint  quadrature rule
\begin{align}\label{eq:semidiscretospaceReaction}
\int_{V_{ij}}  h(u) \dx{x}\dx{y}\approx \lvert V_{ij}\rvert h(\bar{u}_{ij}).
\end{align}

\subsubsection{Implicit part: diffusion space discretization} \label{sec:numSemiDiscImplicit}

The diffusion part is:
\begin{align}
 \oint_{\Gamma_{ij}} {\mbG}(u_x,u_y)\cdot{\mb{ n} } \,d\gamma=& \int_{\Gamma_{i+1/2},j} g_1(u_x,u_y) \dx{\gamma_j} + \int_{\Gamma_{i-1/2,j}} g_1(u_x,u_y) \dx{\gamma_j}  + \notag\\
 & \int_{\Gamma_{i,j+1/2}} g_2(u_x,u_y)  \dx{\gamma_i} + \int_{\Gamma_{i,j-1/2}} g_2(u_x,u_y)  \dx{\gamma_i}. \label{eq:implicitPartDiffusion}
\end{align} 
Before approximating each one of these line integrals, a suitable approximation of the partial derivatives $u_x$ and $u_y$ has to be built for the volume $V_{ij}$. With this aim, we build the second order Lagrange interpolating polynomial of $u$ centered in the volume $V_{ij}$. Let $L_{ij}^u$ be that polynomial. Considering the  nine nodes of the dual mesh $\{x_{i+k},y_{j+l}\},\, k,l=-1,0,1$,  and the averaged values $\{\bar{u}_{i+k,j+l}\}$ of the solution at each volume ${V_{i+kj+l}}$, this polynomial is given by:
$$
L_{ij}^{\bar{u}}(x,y)=\sum_{k,l=-1}^1 \bar{u}_{i+k,j+l} \, \ell_{i+k}(x)\ell_{j+l}(y), 
$$
where $\ell_{i+k}$ and $\ell_{j+l}$ are the one dimensional Lagrange polynomial basis, i.e:
$$
 \ell_{i+k}(x) = \prod_{\substack{p=i-1 \\ p\neq i+k}}^{i+1} \frac{x-x_p}{x_{i+k}-x_p}, \qquad \ell_{j+l}(y) = \prod_{\substack{q=j-1 \\ q\neq j+l}}^{j+1} \frac{y-y_q}{y_{j+l}-y_q}.
$$
Therefore, we use the approximations  $u_x\approx \partial_x L_{ij}^{\bar{u}}$ and $u_y\approx \partial_y  L_{ij}^{\bar{u}}$. The computation of these approximations of the partial derivatives of the solution is straightforward, since it just involves the computation of derivatives of one dimensional Lagrange polynomial basis.

As a summary, the line integrals in equation \eqref{eq:implicitPartDiffusion} are approximated as
\begin{align*}
   \int_{\Gamma_{i\pm 1/2,j}}  & g_1( u_x, u_y)  \dx{\gamma_j}  \approx  \int_{\Gamma_{i\pm 1/2,j}}  g_1\Bigl( \partial_x L_{ij}^{\bar{u}}, \partial_y L_{ij}^{\bar{u}}\Bigr)  \dx{\gamma_j}\approx
\\
& \pm \Delta y\,  g_1\Bigl( \partial_x L_{ij}^{\bar{u}}(x_{i\pm 1/2},y_j  ), \partial_y L_{ij}^{\bar{u}}(x_{i\pm 1/2},y_j  )\Bigr), 
\\
\\
   \int_{\Gamma_{i,j\pm 1/2}}  & g_2( u_x, u_y)  \dx{\gamma_i}  \approx  \int_{\Gamma_{i,j\pm 1/2}}  g_2\Bigl( \partial_x L_{ij}^{\bar{u}}, \partial_y L_{ij}^{\bar{u}}\Bigr)  \dx{\gamma_i}\approx
\\
& \pm \Delta x\,  g_2\Bigl( \partial_x L_{ij}^{\bar{u}}(x_i, y_{j\pm 1/2}  ), \partial_y L_{ij}^{\bar{u}}(x_i, y_{j\pm 1/2} )\Bigr).
\end{align*}

\subsection{IMEX Runge-Kutta time discretization} \label{sec:timeDiscretization}

An IMEX Runge-Kutta scheme consists of applying an implicit discretization to the diffusion terms (stiff terms) and an explicit one to the convective and source terms (non stiff terms). When applied to system \eqref{eq:stiffedo} at second order it takes the form
\begin{align}
U^{(k)} &= U^n -\Delta t \sum_{l=1}^{k-1} \tilde{a}_{kl} {E}(U^{(l)})+\Delta t \sum_{l=1}^{2} a_{kl} {I}(U^{(l)}), \label{eq:RKIMEX_steps}\\
U^{n+1} &= U^n- \Delta t \sum_{k=1}^{2} \tilde{\omega}_{k} E(U^{(k)})+\Delta t \sum_{k=1}^{2} \omega_{k}I(U^{(k)}), \label{eq:RKIMEX_final}
\end{align}
where $U^n=(\bar{u}_{ij}^n)$, $U^{n+1}=(\bar{u}_{ij}^{n+1})$ are the vectors of unknowns volume averages at  times $t^n$ and $t^{n+1}$, and $U^{(k)}$ and $U^{(l)}$ are the vector of unknowns at the stages $k,l$ of the IMEX Runge-Kutta scheme. The matrices $\bar{A} = (\tilde{a}_{ij})$, $\tilde{a}_{kl}=0$ for $l\geq k$ and $A = (a_{kl})$ are $2\times 2$ matrices such that the scheme is explicit in $E$ and implicit in $I$. The coefficient vectors $\tilde{w} = (\tilde{w}_1,\tilde{w}_2)$ and $w = (w_1,w_2)$ complete the IMEX Runge-Kutta scheme.

One should consider diagonally implicit Runge-Kutta (DIRK) schemes \cite{hairerWanner} for the implicit terms ($a_{kl} = 0$, for $l > k$) in order to solve efficiently the algebraic equations corresponding to the implicit part of the discretization at each time step. Besides, the use of a DIRK scheme for the implicit part $I$ guarantees that $F$ is always evaluated explicitly.

IMEX Runge-Kutta schemes can be represented by a double \textit{tableau} in the usual Butcher
notation,

\begin{center}
\begin{tabular}{c|c}
$\tilde{c}$  & $\tilde{A}$ \\
\hline
&  $\tilde{\omega}$ 
\end{tabular}
\quad
\begin{tabular}{c|c}
$c$  & $A$ \\
\hline
         &  $\omega$ 
\end{tabular},
\end{center}
where the coefficients vectors $\tilde{c}=(\tilde{c}_1,\tilde{c}_2)^\mathsf{T}$ and $c=(c_1,c_2)^\mathsf{T}$ used for the treatment of non autonomous systems, are given by the
usual relation
\begin{equation}
  \tilde{c}_k=  \sum_{l=1}^{k-1} \tilde{a}_{kl}, \quad c_k=\sum_{l=1}^{k} {a}_{kl}.
\end{equation}
The PDEs \eqref{eq:PdeBasket2DForward} and \eqref{eq:HestonPDEForward} we will be dealing with in this article, after the space discretization, give rise to a system of non-autonomous ODEs. Therefore, the constants $\tilde{c}_k$ and $c_k$ were omitted in the time marching scheme \eqref{eq:RKIMEX_steps}-\eqref{eq:RKIMEX_final}.

In this work we will consider the second order IMEX-SSP2(2,2,2) (L-stable scheme, see \cite{Russo05}), whose \textit{tableaus} for the explicit (left) and implicit (right) parts are
\begin{center}
\begin{tabular}{c|cc}
0  & 0 & 0 \\
1  & 1 & 0 \\
\hline
& $1/2$ & $1/2$  
\end{tabular}
\quad
\begin{tabular}{c|cc}
$\gamma$  & $\gamma$ & 0 \\
$1-\gamma$  & $1-2\gamma$ & $\gamma$ \\
\hline
& $1/2$ & $1/2$  
\end{tabular}
\quad 
$\gamma=1-\dfrac{1}{\sqrt{2}}.$ 
\end{center}

Typically, an explicit time integrator needs extremely small time steps due to quadratic CFL-like condition:
\begin{align}
    & 2\eta_1 \dfrac{\Delta t}{(\Delta x)^2}  
+2\eta_2 \dfrac{\Delta t}{(\Delta y)^2} 
+\dfrac{1}{2} \eta_3   \dfrac{\Delta t}{\Delta x \Delta y}  
\leq 1, \label{eq:stabilityDiffusion} \\
 & \alpha_1 \dfrac{\Delta t}{\Delta x} + \alpha_2 \dfrac{\Delta t}{\Delta y} \leq 1,\label{eq:stabilityConvection}
\end{align}
where 
\begin{align*}
& \eta_1=\left\lvert\dfrac{\partial g_1}{\partial u_x}\right\rvert,\, \eta_2=\left\lvert\dfrac{\partial g_2}{\partial u_y}\right\rvert,\, \eta_3=\left\lvert\dfrac{\partial g_1}{\partial u_y}+
\dfrac{\partial g_2}{\partial u_x}\right\rvert, \\
& \alpha_1 = \left\lvert\frac{\partial f_1}{\partial u} \right\rvert, \alpha_2 = \left\lvert\frac{\partial f_2}{\partial u} \right\rvert,
\end{align*}
for all volumes $V_{ij}$ and $\forall (x,y)\in\Gamma_{ij}$.
However, IMEX only needs to satisfy the advection stability condition \eqref{eq:stabilityConvection}.

\section{Reference solutions. Fourier Cosine Method\label{2dcos}}

To the best of our knowledge, the most accurate alternative method to assess on the accuracy of our developed finite volume scheme is the COS method proposed by Fang and Oosterlee in \cite{fang:oosterlee:08}. This pricing method is based on Fourier-cosine series expansions. The keystone of Fourier methods is the fact that while probability densities are usually not known for many relevant asset processes, their Fourier transforms, known as the characteristic functions, are available in many relevant cases. Let us describe the method in the two dimensional case.

In this section $t\in[0,T)$ denotes the backward time, so we are interested in the price of the options before the maturity $T$ of the option. Let $\mbx = (x_1,x_2)\in \R^2 $.  Consider a two-dimensional pricing model whose risky assets $\mb{s} = (s_{1},s_{2})$ under a selected martingale measure $Q$ take the form $$\mb{s}_T = (s_{1T},s_{2_T}) = (\e{x_{1T}+x_{10}},\e{x_{2T}+x_{20}})= (s_{10}\e{x_{1T}},s_{20}\e{x_{2T}}),$$ where for $i=1,2$, $x_{i0} = \ln(s_{i0})$, $s_{i0}>0$, and $\mb{x}_T=(x_{1T},x_{2T})$ is a random variable whose characteristic function $\varphi_{\mb{x}_T}$ is known. Let $\mbx_0 = (x_{10},x_{20})$.
The value at time $t$ of an European option with maturity $T$ and payoff function $u_T$ is given by
$$u(\mb{s}_0,t) =   \e{-r(T-t)} \Esp{u_T(\ln(\mb{s}_T))} 
 =  \e{-r(T-t)}\int_{\R^2}u_T(\mbx + \mbx_0)\varPhi_{\mb{x}_T}(\mbx) \dx{\mbx},
$$
where $r$ is the interest rate and $\varPhi_{\mb{x}_T}$ is the probability density function of $\mb{x}_T=(x_{1T},x_{2T})$. Since the tales of $\varPhi_{\mb{x}_T}$  rapidly decay to zero, we can truncate the infinite integration range $\R^2$ to some domain $I=[a_1,b_1]\times[a_2,b_2]\subset \R^2$ without loosing significant accuracy (see \cite{fang:oosterlee:08,ruijter:oosterlee:12} for details on how to reduce the computational domain), thus
$$u(\mb{s}_0,t)  \approx  \e{-r(T-t)}\int_{I} u_T(\mbx + \mbx_0)\varPhi_{\mb{x}_T}(\mbx) \dx{\mbx}.$$
At this point we use the Fourier-cosine approximation of $\varPhi_{\mb{x}_T}$. Let $\mbk = (k_1,k_2)\in\N^2$, and for $i=1,2$ let 
\begin{equation*}
\phi^{(i)}_{k_i}(x_i)= \left\{ \begin{array}{lcc}
             1, &   \mbox{if}  & k_i = 0, \\
             \sqrt{2} \cos\left(k_i \pi \dfrac{x_i-a_i}{b_i-a_i}\right), &  \mbox{if}  & k_i > 0.
             \end{array}
   \right. 
\end{equation*}
The cosine basis functions $\phi_{\mb{k}}(\mb{x}) = \phi^{(1)}_{k_1}(x_1)\phi^{(2)}_{k_2}(x_2)$, are orthonormal in $I$ with respect to the probability density weight $\nu$ of the continuous uniform distribution in $I$, i.e $\nu(\mbx) = \nu_1(x_1)\nu_2(x_2)$, where for $i=1,2$
\begin{equation*} 
\nu_i(x_i)= \left\{ \begin{array}{lcc}
             \dfrac{1}{b_i-a_i}, &   \mbox{if}  & a_i\leq x_i \leq b_i, \\
             \\ 0, &  \mbox{if}  & x_i<a_i \mbox{ or } x_i>b_i.
             \end{array}
   \right. 
\end{equation*}
If $\varPhi_{\mb{x}_T}\in L^2_{\nu}(I)$, it can be decomposed in its Fourier-cosine series
$$\varPhi_{\mb{x}_T}(\mbx) = \sum_{\mbk\in\N^2} \alpha_{\mbk} \phi_{\mbk}(\mbx),$$ where the Fourier coefficients $\alpha_{\mbk}$ are given by $$\alpha_{\mbk} = \int_{I} \varPhi_{\mb{x}_T}(\mbx) \phi_{\mbk}(\mbx) \nu(\mbx) \dx{\mbx}=\left(\prod_{i=1}^2\frac{1}{b_i-a_i} \right)\int_{I} \varPhi_{\mb{x}_T}(\mbx) \phi_{\mbk}(\mbx) \dx{\mbx}.$$ Let $A_{\mbk} = \int_{I} \varPhi_{\mb{x}_T}(\mbx) \phi_{\mbk}(\mbx) \dx{\mbx}$.
In practice, the Fourier coefficients are computed only on a finite set of multi-indices $\Gamma \subset \N^d$, therefore $\varPhi_{\mb{x}_T}(\mbx) \approx \sum_{\mb{k}\in\Gamma} \alpha_{\mbk} \phi_{\mbk}(\mbx)$.
Coming back to the option price, we have that
\begin{align*}
  u(\mb{s}_0,t) &  \approx \e{-r(T-t)} \sum_{\mbk \in \Gamma} \alpha_{\mbk} \int_{I}u_T(\mbx+\mbx_0) \phi_{\mbk}(\mbx)\dx{\mbx} \\
  & = \e{-r(T-t)} \sum_{\mbk \in \Gamma} A_{\mbk} \int_{I}u_T(\mbx+\mbx_0) \phi_{\mbk}(\mbx) \nu(\mbx)\dx{\mbx}.
\end{align*}
Note that $$\beta_{\mbk}(\mbx_0) = \int_{I}u_T(\mbx+\mbx_0) \phi_{\mbk}(\mbx) \nu(\mbx)\dx{\mbx} = \left(\prod_{i=1}^2\frac{1}{b_i-a_i} \right)\int_{I}u_T(\mbx+\mbx_0) \phi_{\mbk}(\mbx) \dx{\mbx},$$ are the Fourier coefficients of the Fourier-cosine expansion of the payoff function $u_T$. Let $B_{\mbk}(\mbx_0) = \int_{I}u_T(\mbx+\mbx_0) \phi_{\mbk}(\mbx) \dx{\mbx}$. All in all we have
$$u(\mb{s}_0,t) \approx \e{-r(T-t)} \left(\prod_{i=1}^2\frac{1}{b_i-a_i} \right) \sum_{\mbk \in \Gamma} A_{\mbk} B_{\mbk}(\mbx_0). $$
The computation of the coefficients $B_{\mbk}(\mbx_0)$ depends on the payoff function $u_T$ and only for a few of these financial functions the integral can be made analytically. For the here discussed arithmetic basket options with payoff function \eqref{eq:arithmeticPayoff}, where $$u_T(x_1,x_2)=\max\left(\frac{1}{2} (s_{10},s_{20})\cdot(\e{x_1},\e{x_2}) ,0\right),$$ the payoff Fourier coefficients have to be computed numerically. Let us now  discuss how to compute the coefficients $A_{\mbk}$ efficiently. For the sake of brevity we restrict the notation to the two dimensional case. In order to alleviate the notation, let us denote $\hat{k}_i = \frac{\pi}{b_i-a_i}k_i \in \R$ for $i=1,2$. In this case,
\begin{equation*} 
\phi_{\mbk}(\mbx)= \left\{ \begin{array}{lcl}
             1, &   \mbox{if}  & \mbk=(0,0), \\
             \sqrt{2} \cos\left(\hat{k}_1 (x_1-a_1)\right), &  \mbox{if}  & \mbk=(k_1,0),\, k_1>0,\\
             \sqrt{2} \cos\left(\hat{k}_2 (x_2-a_2)\right), &  \mbox{if}  & \mbk=(0,k_2), \, k_2>0,\\
             2 \cos\left(\hat{k}_1 (x_1-a_1)\right) \cos\left(\hat{k}_2 (x_2-a_2)\right), &  \mbox{if}  & \mbk=(k_1,k_2),\,k_1,k_2>0.                         
             \end{array}
   \right. 
\end{equation*}
Since
\begin{align*}
 &\cos\left(\hat{k}_i (x_i-a_i)\right) = \Real{\e{-\mb{i}a_i \hat{k}_i} \e{\mb{i}x_i \hat{k}_i}}, \mbox{ and}\\
 &2 \cos\left(\hat{k}_1 (x_1-a_1)\right) \cos\left(\hat{k}_2 (x_2-a_2)\right) = \\
 &\,\,  \cos\left(\hat{k}_1 (x_1-a_1) + \hat{k}_2 (x_2-a_2)\right) + \cos\left(\hat{k}_1 (x_1-a_1) - \hat{k}_2 (x_2-a_2)\right) = \\
 & \,\, \Real{\e{-\mbi(a_1,a_2)\cdot\left( \hat{k}_1,\hat{k}_2 \right) } \e{\mbi(x_1,x_2)\cdot\left( \hat{k}_1,\hat{k}_2 \right) }  } + \Real{\e{\mbi(a_1,a_2)\cdot\left( -\hat{k}_1,\hat{k}_2 \right) } \e{\mbi(x_1,x_2)\cdot\left( \hat{k}_1,-\hat{k}_2 \right) }  }, \\
 & \mbox{and} \\
 & \int_{I} \varPhi_{\mb{x}_T}(\mbx) \e{\mbi (x_1,x_2)\cdot(\hat{k}_1,\hat{k}_2) } \dx{\mbx} \approx \int_{\R^2} \varPhi_{\mb{x}_T}(\mbx) \e{\mbi (x_1,x_2)\cdot(\hat{k}_1,\hat{k}_1) } \dx{\mbx} = \varphi_{\mb{x}_T}\left(\hat{k}_1,\hat{k}_1\right),
\end{align*}
it is clear that 
\begin{equation} \label{eq:aCoeffs}
A_{(k_1,k_2)}\approx \left\{ \begin{array}{lcl}
             \Real{\varphi_{\mb{x}_T}(0,0)}, &   \mbox{if}  & \mbk=(0,0), \\
             \sqrt{2}\Real{\e{-\mbi a_1\hat{k}_1} \varphi_{\mb{x}_T}\left(\hat{k}_1,0\right) }, &  \mbox{if}  & \mbk=(k_1,0),\, k_1>0,\\
             \sqrt{2}\Real{\e{-\mbi a_2\hat{k}_2} \varphi_{\mb{x}_T}\left(0, \hat{k}_2 \right) }, &  \mbox{if}  & \mbk=(0,k_2), \, k_2>0,\\
             \Real{ \e{-\mbi(a_1,a_2)\cdot\left( \hat{k}_1,\hat{k}_2 \right) }  \varphi_{\mb{x}_T}\left( \hat{k}_1,\hat{k}_2 \right) } + \\
             \,\, \Real{ \e{\mbi(a_1,a_2)\cdot\left( -\hat{k}_1,\hat{k}_2 \right) }  \varphi_{\mb{x}_T}\left( \hat{k}_1,-\hat{k}_2 \right) } , &  \mbox{if}  & \mbk=(k_1,k_2),\,k_1,k_2>0.
             \end{array}
   \right. 
\end{equation}

This spectral method suffers the curse of dimensionality: not only in the computation of the multidimensional integrals defining $B_{\mbk}(\mbx_0)$, but also in the series expansion itself, especially if $\varPhi_{\mb{x}_T}$ is a heavy tailed density. A smart choice of the multi-indices coefficients $\Gamma$ will help to mitigate the curse of dimensionality of this algorithm.

In the Black-Scholes model \eqref{eq:sdesBS2d} we have
\begin{equation} \label{eq:characteristicBS}
 \varphi_{\mb{X}_T}(\mb{u}) = \exp\left(  \mbi T \boldsymbol{\mu} \mb{u}^\mathsf{T} - \frac{1}{2} T \mb{u}\Sigma \mb{u}^\mathsf{T} \right),
\end{equation}
where $\boldsymbol{\mu}=(\mu_1,\mu_2)$ with $\mu_i = (r-q_i-\frac{1}{2}\sigma_i^2)$ and $\Sigma$ is the covariance matrix with elements $\Sigma_{i,j}=\sigma_i \sigma_j \rho_{ij}$. Let $c^{(i)}_1 = T \left( r-\frac{1}{2}\sigma_i^2\right)$ and $c^{(i)}_2 = T \sigma_i^2$ be the cumulants of the Brownian motion. Here, the way to reduce the computational domain is to consider (see \cite{fang:oosterlee:08})
\begin{equation} \label{eq:cutCompDom}
a_i = c^{(i)}_1 - L\sqrt{c^{(i)}_2},\quad b_i = c^{(i)}_1 + L\sqrt{c^{(i)}_2},\quad L\in\R_+. 
\end{equation}

At this point let us discuss the computation of the coefficients $A_{\mbk}$ in the general case ($k_1,k_2>0$) for the Black-Scholes model. By plugging \eqref{eq:characteristicBS} into \eqref{eq:aCoeffs} considering \eqref{eq:cutCompDom}, since $\frac{T\mu_i-a_i}{b_i-a_i}=\frac{1}{2} $, one readily obtains
\begin{align*}
  & A_{(k_1,k_2)} \approx  \\
 & \, \Real{ \e{-\mbi(a_1,a_2)\cdot\left(\hat{k}_1,\hat{k}_2  \right) }  \varphi_{\mb{x}_T}\left( \hat{k}_1,\hat{k}_2 \right) } + \Real{ \e{\mbi(a_1,a_2)\cdot\left( -\hat{k}_1,\hat{k}_2 \right) }  \varphi_{\mb{x}_T}\left( \hat{k}_1,-\hat{k}_2 \right) }\\
 & \,=  C_1 \Real{\e{\mbi \left(T\mu_1-a_1,T\mu_2-a_2 \right) \left(\hat{k}_1,\hat{k}_2\right) }} +  C_2 \Real{\e{\mbi \left(T\mu_1-a_1,T\mu_2-a_2 \right) \left(\hat{k}_1,-\hat{k}_2\right) }} \\
 & \,=  C_1 \cos\left( (k_1+k_2)\frac{\pi}{2} \right) + C_2 \cos\left( (k_1-k_2)\frac{\pi}{2} \right),
\end{align*}
where $C_1$, $C_2 \in \R$ are given by
\begin{align*}
 C_1 &= \exp\left(  - \frac{1}{2} T \left(\hat{k}_1,\hat{k}_2\right) \Sigma \left(\hat{k}_1,\hat{k}_2\right)^\mathsf{T}\right), \\
 C_2 &= \exp\left(  - \frac{1}{2} T \left(\hat{k}_1,-\hat{k}_2\right) \Sigma \left(\hat{k}_1,-\hat{k}_2\right)^\mathsf{T}\right).
\end{align*}
As a consequence, for the Black-Scholes model $A_{(k_1,k_2)}\neq 0$ if and only if $k_1+k_2$ is even. Therefore, the here considered set of multi-indices is $$\Gamma^{BS} = \left\{ (k_1,k_2)\in\N^2: k_1\leq K_1, k_2\leq K_2, k_1+k_2 \mbox{ even} \right\}.$$ This choice will allow us to overcome the curse of dimensionality in the series expansion. Thus, we will be able to compute, at the same cost, more accurate option prices. Regarding the numerical computation of the payoff coefficients $B_{\mbk}(\mbx_0)$, we have used composed Gauss-Legendre quadrature rule with two thousand points on each spatial direction. In order to speedup the computations, we have taken advantage of high performance computing techniques. The loop for the multi-indices ${\mbk}\in\Gamma$ was parallelized in CPU using four OpenMP threads. Each one of those threads, is responsible for computing $A_{\mbk}$ and $B_{\mbk}(\mbx_0)$. The calculation of $B_{\mbk}(\mbx_0)$ is very costly. In order to accelerate this computation, we have parallelized the quadrature rule using Nvidia CUDA code. This code was ran on four Nvidia Tesla V100.

\algblock{ParFor}{EndParFor}
\algnewcommand\algorithmicparfor{\textbf{parfor}}
\algnewcommand\algorithmicpardo{\textbf{do}}
\algnewcommand\algorithmicendparfor{\textbf{end\ parfor}}
\algrenewtext{ParFor}[1]{\algorithmicparfor\ #1\ \algorithmicpardo}
\algrenewtext{EndParFor}{\algorithmicendparfor}

\begin{algorithm}
\caption{Sketch of parallel COS method}\label{algo1}
\begin{algorithmic}[1]
\Require $s_0 > 0$
\State $u(\mb{s}_0,t) = 0$
\ParFor{$\mbk\in\Gamma^{BS}$} \algorithmiccomment{This is a parallel OpenMP for loop}
\State Compute $A_{\mbk}$
\State Compute $B_{\mbk}(\mbx_0)$  \algorithmiccomment{GPU kernel in NVIDIA V100}
\State $u(\mb{s}_0,t)  += A_{\mbk} B_{\mbk}(\mbx_0)$ \algorithmiccomment{Reduction variable}
\EndParFor
\State $u(\mb{s}_0,t) = \e{-r(T-t)} \left(\prod_{i=1}^2\frac{1}{b_i-a_i} \right) u(\mb{s}_0,t)$
\end{algorithmic}
\end{algorithm}

Finally, in order to price options under the Heston model using the COS method, several approaches are available. For example, one could use the two dimensional COS here presented, using the bivariate characteristic function presented in \cite{ruijter:oosterlee:12}. Nevertheless, in the Fourier framework, Heston problem can be reduced to a one-dimensional problem. Here we have considered the following characteristic function 
\begin{align*}
 \varphi_{x_T}(\xi) =& \exp\left( \mbi \xi (r-q) T+ \frac{\nu_0}{\sigma^2} \left( \frac{1-\e{-D(\xi)T}}{1-G(\xi)\e{-D(\xi)T}} \right)  (\kappa - \mbi \rho \sigma \xi - D(\xi))\right) \times \\
 & \exp\left( \frac{\kappa \theta}{\sigma^2} \left( T (\kappa - \mbi \rho \sigma \xi - D(\xi)) - 2 \ln\left(\frac{1-G(\xi)\e{-D(\xi)T}}{1-G(\xi)} \right)\right)\right),
\end{align*}
where $\nu_0$ is the initial variance and $D$, $G$ the following functions
\begin{equation*}
 D(\xi) = \sqrt{(\kappa-\mbi \rho \sigma \xi)^2 + (\xi+\mbi)\xi \sigma^2}, \quad G(\xi) = \frac{\kappa-\mbi \rho \sigma \xi - D(\xi)}{\kappa - \mbi \rho \sigma \xi + D(\xi)}.
\end{equation*}
Besides, the payoff coefficients for the vanilla call option with payoff $\max(s_T-K,0)$ can be computed analytically (see \cite{fang:oosterlee:08} for details).

\section{Numerical experiments \label{sec:numerical-results}}

In this section, several numerical experiments are developed to assess the accuracy and performance of the new numerical methods proposed in previous sections for the discussed two dimensional problems. 

Experiments are carried out on several option pricing problems. In subsection \ref{sec:NumExpBasket} basket call options over two underlyings following Black-Scholes model are priced. Vanilla call options under Heston stochastic volatility model are priced in subsection \ref{sec:NumExpHeston}. Each subsection is organized as follows. 

We start by writing the PDE model in conservative form. Then, we solve the PDE problems using the discussed finite volume IMEX Runge-Kutta numerical method. Both convection-dominated and diffusion-dominated scenarios are considered. As a result, surface prices are presented. In addition, plots for the numerical derivatives of the solution are presented. More precisely, delta and gamma Greeks are shown. All pictures, tables and errors in this whole section will be computed considering the numerical solution at the last time step where $t=T$, i.e, we approximate numerically the option prices ``today''.

Later, reference solutions for the studied pricing problems are accurately computed at $t=T$ by means of the COS Fourier method explained in Section \ref{2dcos}. They allow us to compute $L_1$ and $L_\infty$ error norms, and the following relative and mean absolute errors given by:
\begin{align*}
    L_1 \mbox{ error } &= \Delta x \Delta y\sum_{i,j}  \lvert \bar{u}_{i,j}-u(x_i,y_j) \rvert, \\
    L_\infty \mbox{ error } &= \max_{i,j} \lvert \bar{u}_{i,j}-u(x_i,y_j) \rvert, \\    
    L_\infty \mbox{ relative error } &= \dfrac{L_\infty \mbox{ error }}{\max_{i,j} \lvert u(x_i,y_j) \rvert},\\
    \mbox{Mean absolute error} &=\underset{i, j}{ \mbox{mean} } \lvert \bar{u}_{i,j}-u(x_i,y_j) \rvert.
\end{align*}
Besides, absolute errors of the finite volume IMEX Runge-Kutta method against the COS Fourier method are plotted in order to quickly comprehend the quality of the approximation and the regions with higher errors. 

On top of that, orders of convergence are validated in-depth. Both explicit and IMEX Runge-Kutta time integrators were implemented over the same explained space discretization. In the explicit case second order Heun's ODE solver was considered. Finite volume with both time marching schemes achieves the announced second order accuracy in all regimes. IMEX and explicit methods yield similar results in terms of accuracy and convergence order. Nevertheless, the IMEX solver can take advantage of much larger times steps, thus offering much better performance in terms of execution times. The explicit method dramatically slows down the computation due to the presence of diffusive terms.

Finally, each subsection is closed with a table with some reference option prices. Since the here considered option prices are not known is closed form, an interested reader may find this table useful, particularly those developers of alternative option price calculators.

Note that the designed algorithms were implemented using C++ (GNU C++ compiler 9.3.0) and were ran in a machine with 128 GB of RAM, one AMD Ryzen 9 5950X 16-Core Processor, and four GPUs NVIDIA Volta V100. All codes were compiled using double precision arithmetic. Besides, for linear algebra operations Eigen library \cite{eigenweb} was considered. For all experiments in this article a CFL of $0.5$ is taken into account in the stability conditions.

\subsection{Options on a basket of two assets\label{sec:NumExpBasket}}

In this experiment we solve the basket European option pricing problem considering two underlyings. The PDE model \eqref{eq:PdeBasket2DForward} can be written in the conservative form \eqref{SistCons} as follows
\begin{equation} \label{SistConsBS}
    \dfrac{\partial u}{\partial t}  + \dfrac{\partial f_1}{\partial s_1} (u)+\dfrac{\partial f_2}{\partial s_2}  (u)=\dfrac{\partial g_1}{\partial s_1} (u_{s_1},u_{s_2})+\dfrac{\partial g_2}{\partial s_2} (u_{s_1},u_{s_2})+h(u),
\end{equation}
where the functions $f_1, f_2, g_1, g_2$ and $h$ are given by:
\begin{align*}
f_1(u)&=(\sigma_1^2-r+q_1)s_1u(s_1,s_2,t)+\frac{\rho}{2}\sigma_1\sigma_2s_1 u(s_1,s_2,t), \\
f_2(u)&=(\sigma_2^2-r+q_2)s_2u(s_1,s_2,t)+\frac{\rho}{2}\sigma_1\sigma_2s_2 u(s_1,s_2,t), \\
g_1(u_{s_1},u_{s_2})&=\frac12\sigma_1^2s_1^2u_{s_1}(s_1,s_2,t)+\frac\rho2\sigma_1\sigma_2s_1s_2 u_{s_2}(s_1,s_2,t),\\
g_2(u_{s_1},u_{s_2})&=\frac12\sigma_2^2s_2^2 u_{s_2}(s_1,s_2,t)+\frac\rho2\sigma_1\sigma_2s_1s_2 u_{s_1}(s_1,s_2,t),\\
h(u)&=( \sigma_2^2+\rho \sigma_1\sigma_2+\sigma_1^2+q_1+q_2-3r ) u(s_1,s_2,t).
\end{align*}
The initial condition $u_0$, being equation \eqref{eq:arithmeticPayoff} for the arithmetic basket call option, was averaged on each volume. Indeed, the initial volume averages are set as
$${\bar{u}^0}_{i,j} = \dfrac{1}{\lvert V_{ij} \rvert} \int_{V_{ij}} u_0(s_1,s_2)\dx{s_2}\dx{s_1},$$
where the mid-point quadrature formula is used to approximate the integral.

In this section we perform two basket option pricing tests with the market parameters of Table \ref{tb-Basket}. The first set of parameters, denoted as Test 1, stands for a convection-dominated regime. Although nowadays this setup with an interest rate $r=0.5$ is financially unrealistic, it is useful as a stress-test. The second group of parameters, labelled as Test 2, represents a diffusion-dominated setting. In both scenarios, in order to avoid noise coming from the artificial boundary conditions, the computational space domain was considered as $(s_1,s_2)\in[0,5K]\times[0,5K]$.

\begin{table}[h]
\begin{center}
\caption{Market data for basket options under Black-Scholes model.}\label{tb-Basket}%
\begin{tabular}{c|c|c|c|c|c|c|c|c}
  &$\sigma_1$ & $\sigma_2$ & $r$  & $q_1$ & $q_2$  & $\rho$ & $T$& $K$\\
\hline
Test 1& $0.1$ & $0.1$  & $0.5$ & \multirow{2}{*}{$0.00$} & \multirow{2}{*}{$0.00$} & \multirow{2}{*}{$0.5$} & \multirow{2}{*}{$0.25$} & \multirow{2}{*}{$30$}\\

Test 2 &$0.5$ & $0.5$  & $0.1$ &  &  &  &  & \\
\end{tabular}
\end{center}
\end{table}

\begin{figure}[!htb]
\centering
\subfigure[Test 1] {\includegraphics[height=4.4cm]{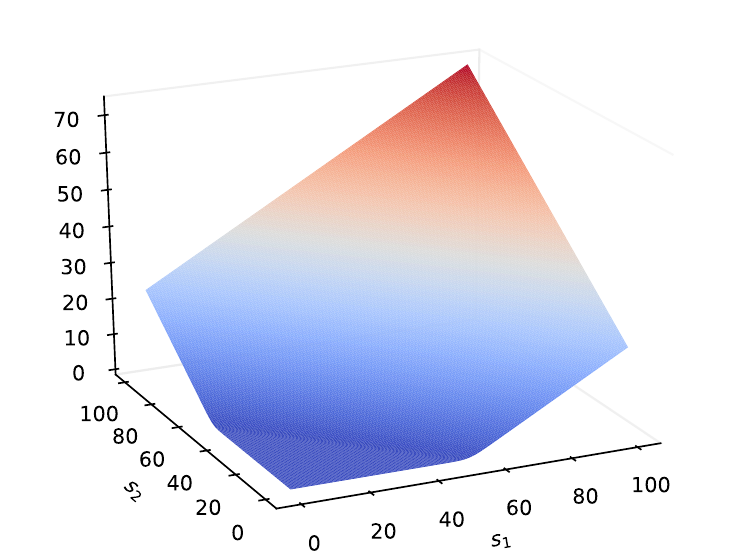}}
\subfigure[Test 2] {\includegraphics[height=4.4cm]{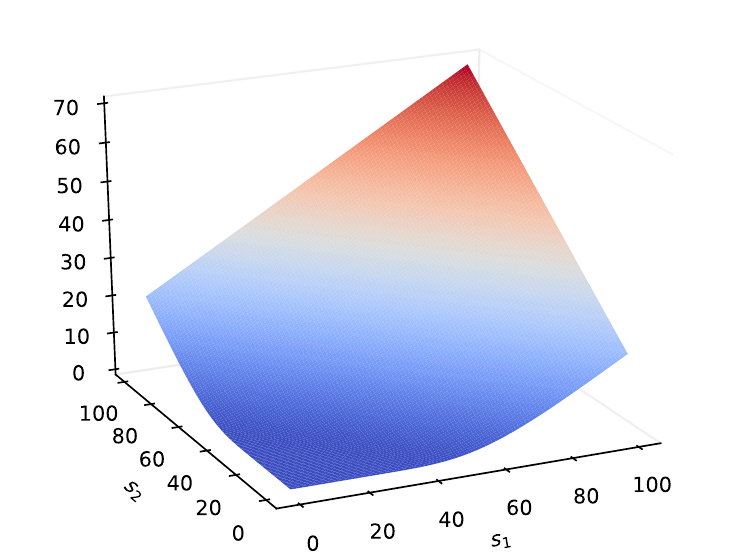}}
\caption{Basket option prices for Test 1 (left) and Test 2 (right) at $t=T$.}
\label{fig:basket-prices}
\end{figure}

\begin{figure}[!htb]
\centering
\subfigure[Test 1] {\includegraphics[height=4.4cm]{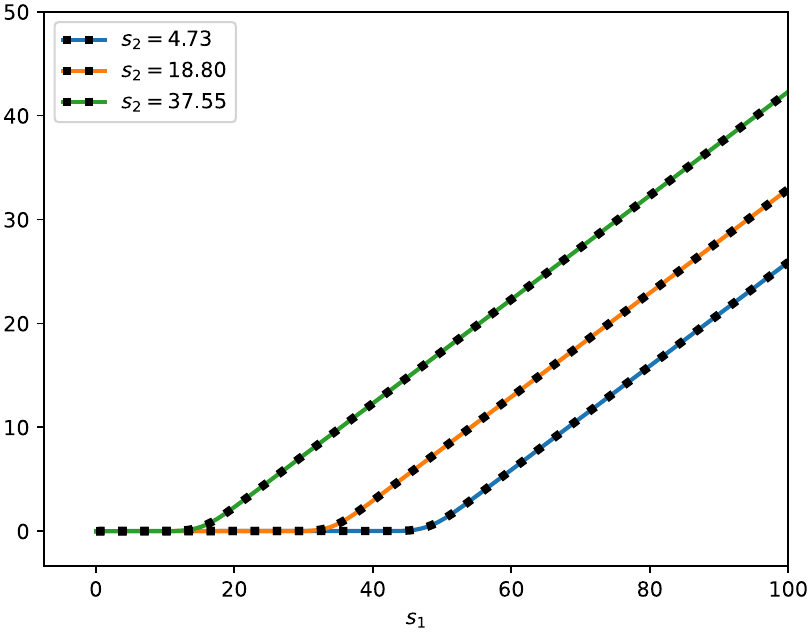}}
\subfigure[Test 2] {\includegraphics[height=4.4cm]{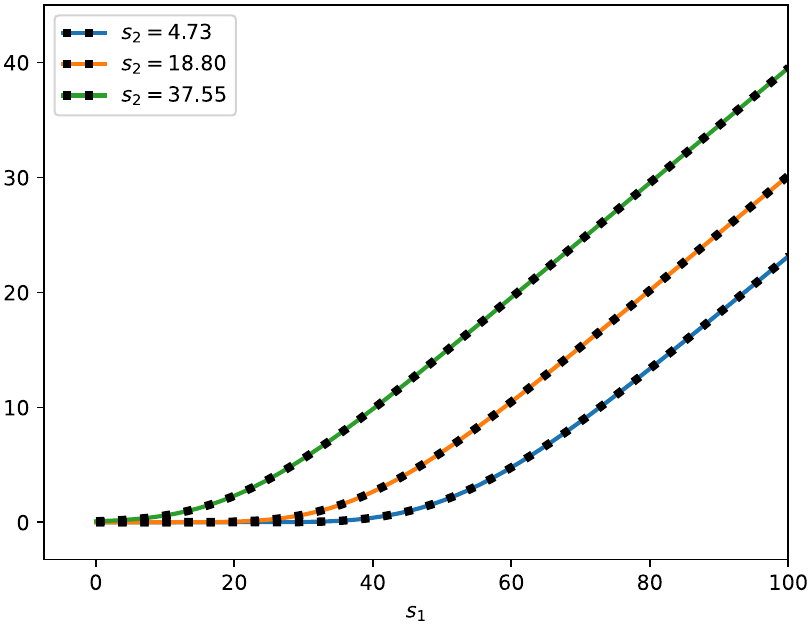}}
\caption{Cuts of price surfaces of Figure \ref{fig:basket-prices}. Numerical solution is shown in continuous line and COS solution with squares.}
\label{fig:basket-prices-cuts}
\end{figure}

\begin{figure}[!htb]
\centering
\subfigure[Test 1] {\includegraphics[height=4.4cm]{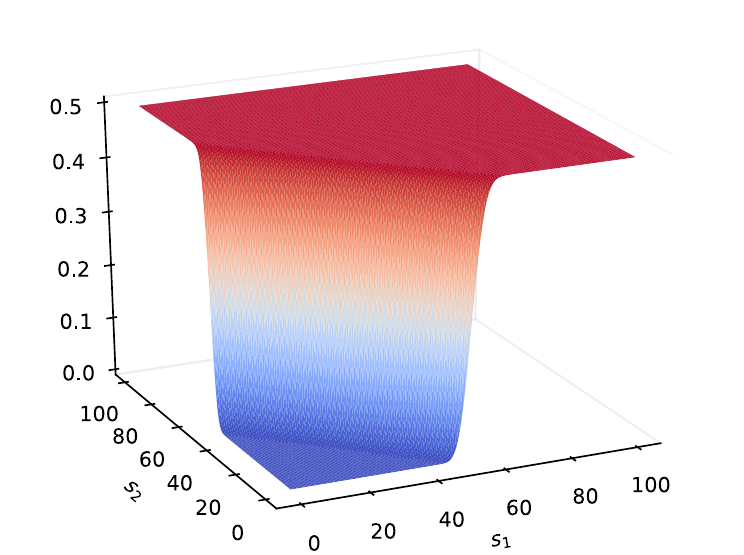}}
\subfigure[Test 2] {\includegraphics[height=4.4cm]{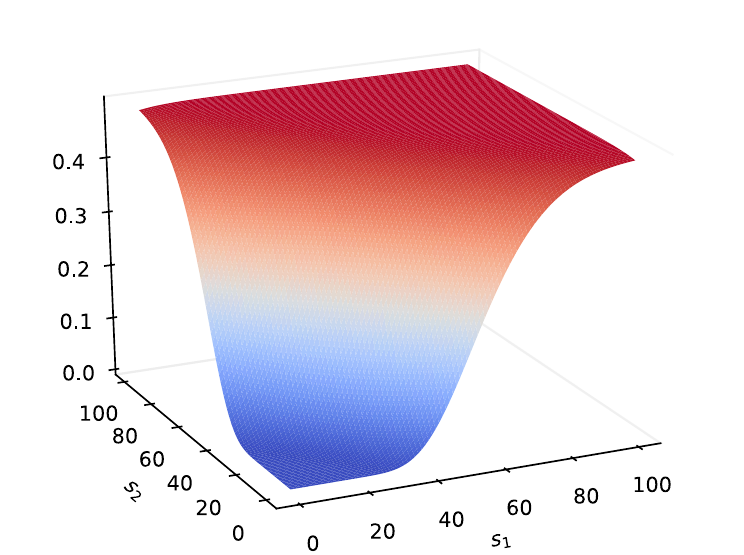}}
\caption{Basket option deltas for Test 1 (left) and Test 2 (right) at $t=T$.}
\label{fig:basket-delta}
\end{figure}

\begin{figure}[!htb]
\centering
\subfigure[Test 1] {\includegraphics[height=4.4cm]{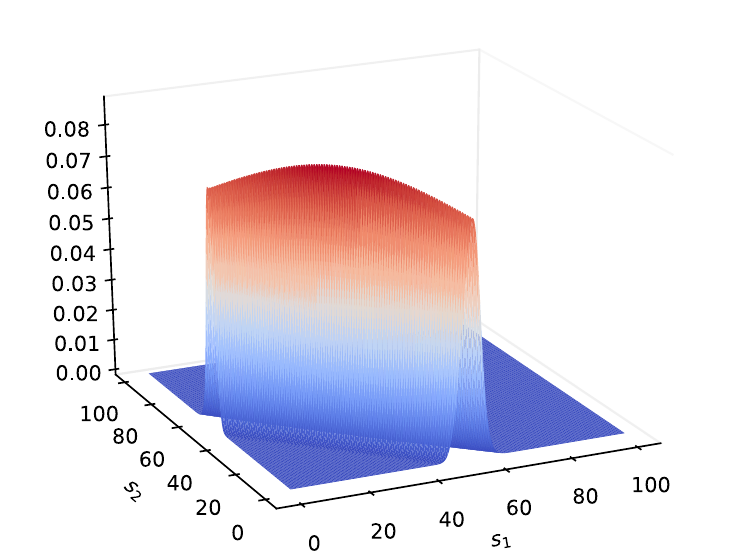}}
\subfigure[Test 2] {\includegraphics[height=4.4cm]{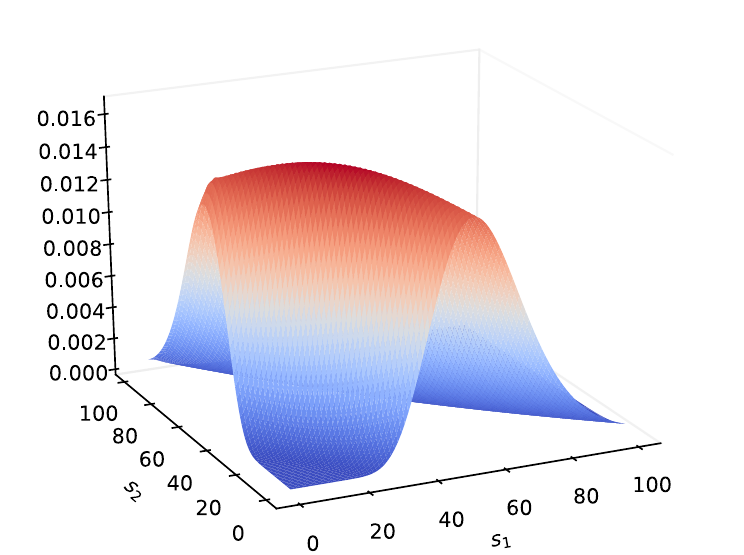}}
\caption{Basket option gammas for Test 1 (left) and Test 2 (right) at $t=T$.}
\label{fig:basket-gamma}
\end{figure}

Surfaces of basket option prices at $t=T$ for Test 1 and Test 2 are shown in Figure \ref{fig:basket-prices}. Cuts of these surfaces with planes parallel to $s_2=0$ are shown in Figure \ref{fig:basket-prices-cuts}, along with COS solutions. All these plots were drawn using the solution obtained with the finite volume IMEX Runge-Kutta solver. Moreover, first numerical derivative (delta) and second numerical derivative (gamma) of the numerical approximation of the option prices are shown in Figures \ref{fig:basket-delta} and \ref{fig:basket-gamma}, respectively. The finite volume IMEX Runge-Kutta numerical scheme offers high resolution approximation of basket option values. Furthermore, it provides accurate and non-oscillatory approximations of the Greeks.

\begin{figure}[!htb]
\centering
\subfigure[Test 1] {\includegraphics[height=4.4cm]{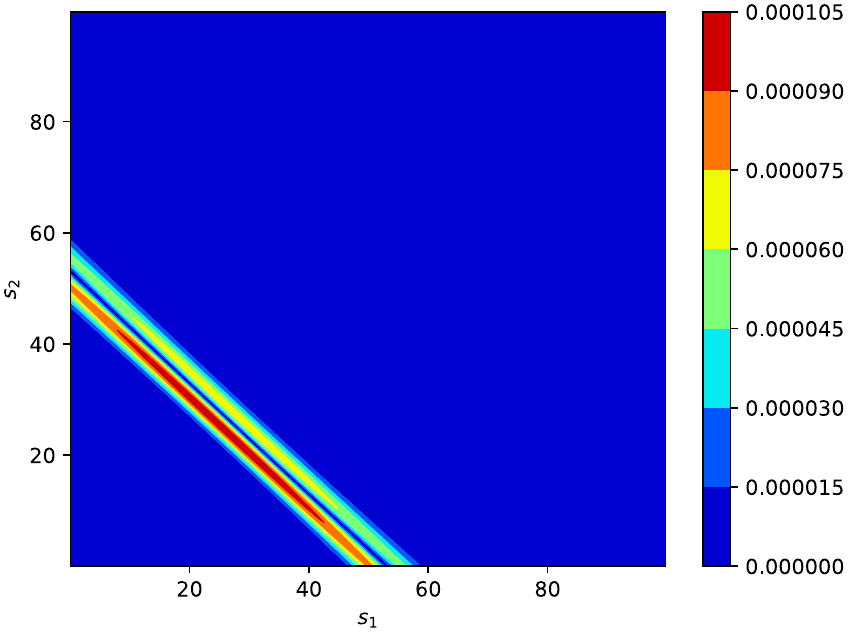}}
\subfigure[Test 2] {\includegraphics[height=4.4cm]{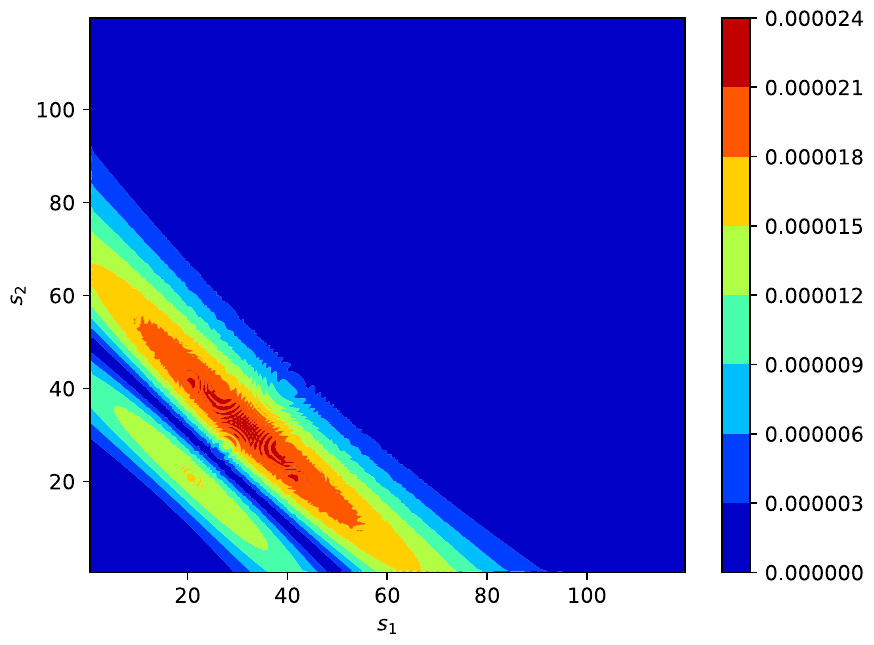}}
\caption{Contour plots for the surface of absolute errors of basket options prices for Test 1 (left) and Test 2 (right) at $t=T$.}
\label{fig:error_basket}
\end{figure}

At this point we compute the reference option prices by using the described COS method. For a mesh of size $N_1\times N_2 = 1600\times1600$, i.e., a mesh with 1600 discretization points in each space direction, the errors of the approximation obtained with the finite volume IMEX Runge-Kutta method are collected in Table \ref{tb-Basket_errors}. Figure \ref{fig:error_basket} presents contour plots for the surface defined by the absolute errors in the option prices at $t=T$. The here developed finite volumes (explicit or IMEX) schemes offer high-resolution approximations even at regions of discontinuities and non-smoothness in the initial condition. Therefore, there is no need to apply smoothing methods like the Rannacher technique.

\begin{table}[!htb]
\begin{footnotesize}
\begin{center}
\caption{Finite volume IMEX Runge-Kutta numerical errors against reference basket option prices computed with COS method at $t=T$.}
\begin{tabular}{|c|c|c|}
\hline
\multicolumn{3}{|c|}{Test 1}\\
\hline
  $L_\infty  \mbox{ error }$ & $L_\infty \mbox{ relative error }$ & $\mbox{Mean absolute error}$ \\
\hline
$1.0245\times 10^{-4}$ & $1.3976\times 10^{-6}$ & $3.5079\times 10^{-6}$  \\
\hline

\hline
\multicolumn{3}{|c|}{Test 2}\\
\hline
  $L_\infty  \mbox{ error }$ & $L_\infty \mbox{ relative error }$ & $\mbox{Mean absolute error}$ \\
\hline
$2.3406\times 10^{-5}$ & $1.9439\times 10^{-7}$ & $2.7099\times 10^{-6}$  \\
\hline
\end{tabular} 
\label{tb-Basket_errors} 
\end{center}
\end{footnotesize}
\end{table}

\begin{figure}[!h]
\centering\includegraphics[width=0.5\linewidth]{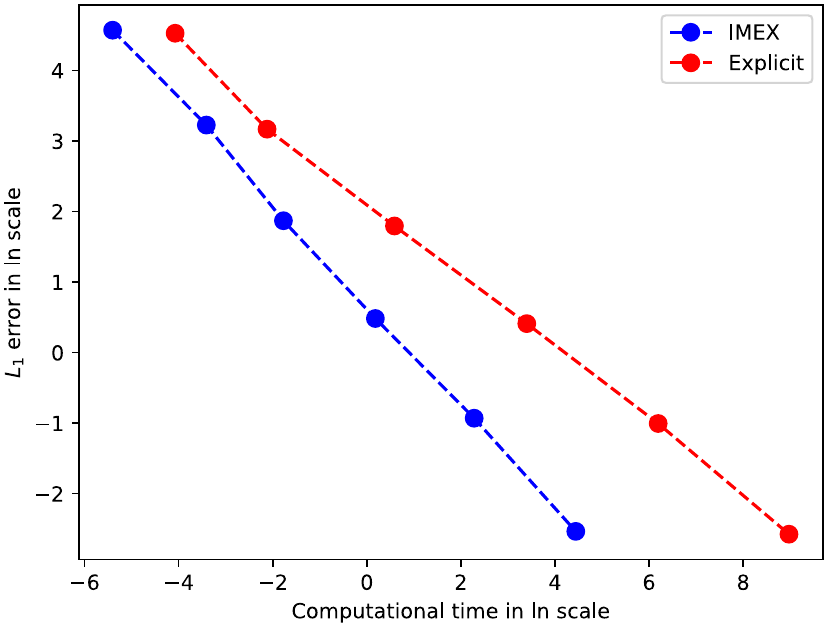}  
\caption{Efficiency curve of IMEX and explicit time marching schemes for Test 1.}
\label{fig:curveEfficiency}
\end{figure}

Tables \ref{tb-Basket_convective} and \ref{tb-Basket_diffusive} record $L_1$ errors and $L_1$ orders of convergence for Test 1 and Test 2, respectively. The errors and orders of convergence are shown for both the IMEX and explicit finite volume numerical schemes. Both schemes achieve second-order accuracy in the $L_1$ norm. Additionally, these tables present the time steps and the execution times for the two time integrators. As said before, IMEX method is able to converge using much larger time steps than the explicit one. The stability condition of the explicit scheme requires extremely small time steps, thus making the method useless in practice for refined grids in space. Therefore, IMEX offers much better performance in terms of execution times, and allows us to solve the PDE problems with refined meshes in space. In fact, IMEX method is able to run for grids finer than $800\times 800$ grid, while the explicit method is not in reasonable computational times. In Figure \ref{fig:curveEfficiency}, the natural logarithms of $L_1$ errors and execution times of Table \ref{tb-Basket_diffusive} are shown for both the IMEX and explicit numerical schemes. IMEX clearly outperforms the explicit method in this scenario with diffusion dominance, which is the typical situation in finance. Although for convective dominated problems in coarse grids both time marching schemes perform similarly, as soon as the mesh is refined IMEX is the only practical choice.

\begin{table}[!h]
\begin{center}
{\scriptsize
\begin{tabular}{|c||c|c|c|c|}
\hline
& \multicolumn{4}{c|}{IMEX}  \\
\hline
$N_1\times N_2$  & $L_1$ \text{ error} &  \text{ Order}& $\Delta t$ & Time (s)  \\
\hline
$25\times 25$ &  $9.9867 \times 10^{1}$ & $--$ & $2.05\times 10^{-2}$ & $1.1\times 10^{-2}$ \\
$50\times 50$ &  $3.3457 \times 10^{1}$ & $1.57$& $1.03\times 10^{-2}$ & $3.9\times 10^{-2} $ \\ 
$100\times 100$ &  $9.1341 \times 10^{0}$ & $1.87$ & $5.13\times 10^{-3}$ & $2.3\times 10^{-1}$ \\ 
$200\times 200$ &  $2.3529 \times 10^{0}$ & $1.95$ & $2.56\times 10^{-3}$ & $1.8\times 10^{0}$ \\
$400\times 400$ &  $5.6234 \times 10^{-1}$ & $2.06$ & $1.28\times 10^{-3}$ & $1.6\times 10^{1}$ \\
$800\times 800$ &  $1.1257\times 10^{-1}$ & $2.32$ & $6.41\times 10^{-4}$ & $1.5\times 10^{2}$ \\ 
\hline
\hline
& \multicolumn{4}{c|}{Explicit}  \\
\hline
$N_1\times N_2$  & $L_1$ \text{ error} &  \text{ Order}& $\Delta t$ & Time (s)  \\
\hline
$25\times 25$ & $9.9752\times 10^{1}$ & $--$ &   $3.56\times 10^{-2}$ & $9.4\times 10^{-4} $\\
$50\times 50$ & $3.3453\times 10^{1}$ &$1.57$& $8.88\times 10^{-3}$ & $1.1\times 10^{-2} $\\
$100\times 100$ & $9.1503 \times 10^{0}$ &$1.87$ & $2.22\times 10^{-3}$ & $8.3\times 10^{-2}$\\
$200\times 200$ & $2.3614\times 10^{0}$ &$1.95$ & $5.56\times 10^{-4}$ &$ 1.4\times 10^{0}$\\
$400\times 400$ & $5.6466\times 10^{-1}$ &$2.06$ &  $1.39\times 10^{-4}$ & $2.1\times 10^{1}$\\
$800\times 800$ & $1.1304\times 10^{-1}$ &$2.32$ &  $3.47\times 10^{-5}$ & $3.6\times 10^{2}$ \\
\hline
\end{tabular} 
}
\end{center}
\caption{$L_1$ errors and $L_1$ orders of convergence of the IMEX and explicit finite volume methods for Test 1.}
\label{tb-Basket_convective} 
\end{table}

\begin{table}[!h]
\begin{center}
{\scriptsize
\begin{tabular}{|c||c|c|c|c|}
\hline
& \multicolumn{4}{c|}{IMEX}  \\
\hline
$N_1\times N_2$  & $L_1$ \text{ error} &  \text{ Order}& $\Delta t$ & Time (s) \\
\hline
$25\times 25$ & $9.6620\times 10^{1}$ & $--$ & $4.71\times 10^{-2}$ & $4.5\times 10^{-3}$ \\
$50\times 50$ &  $2.5178\times 10^{1}$ & $1.94$& $2.35\times 10^{-2}$ & $3.3\times 10^{-2} $ \\
$100\times 100$ &  $6.4828\times 10^{0}$ & $1.95$& $1.18\times 10^{-2}$ & $1.7 \times 10^{-1} $ \\
$200\times 200$ &  $1.6209\times 10^{0}$ & $2.00$ & $5.89\times 10^{-3}$ & $1.2\times 10^{0} $ \\
$400\times 400$ &  $3.9419\times 10^{-1}$ & $2.03$ & $2.94\times 10^{-3}$ & $9.8\times 10^{0}$ \\
$800\times 800$ &  $7.9229\times 10^{-2}$ & $2.31$ & $1.47\times 10^{-3}$ & $8.5\times 10^{1}$ \\
\hline
\hline
& \multicolumn{4}{c|}{Explicit}  \\
\hline
$N_1\times N_2$  & $L_1$ \text{ error} &  \text{ Order}& $\Delta t$ & Time (s) \\
\hline 
$25\times 25$ & $9.0224\times 10^{1}$ & $--$ & $1.42\times 10^{-3}$ & $1.7\times 10^{-2}$\\
$50\times 50$ & $2.3440\times 10^{1}$ &$1.94$ &  $3.56\times 10^{-4}$ & $1.2\times 10^{-1}$\\
$100\times 100$ & $5.9498\times 10^{0}$ &$1.97$ & $8.89\times 10^{-5}$ & $1.8\times 10^{0}$\\
$200\times 200$ & $1.4834\times 10^{0}$ &$2.00$ & $2.22\times 10^{-5}$ & $3.0\times 10^{1}$\\
$400\times 400$ & $3.5473\times 10^{-1}$ &$2.06$ & $5.56\times 10^{-6}$ & $4.9 \times 10^{2}$\\
$800\times 800$ & $7.1095\times 10^{-2}$ &$2.31$ & $1.39\times 10^{-6}$ & $ 7.9 \times 10^{3}$ \\
\hline
\end{tabular} 
}
\end{center}
\caption{$L_1$ errors and $L_1$ orders of convergence of the IMEX and explicit finite volume methods for Test 2.}
\label{tb-Basket_diffusive}
\end{table}

Finally, Tables \ref{tb-basket-convective_anhelada} and \ref{tb-basket-diffusive_anhelada} show some basket option prices at $t=T$, for Test 1 and Test 2, respectively. These prices were obtained with COS method and with the finite volume IMEX Runge-Kutta solver executed over a space grid of size $1600\times 1600$. In the COS method, in order to truncate the series expansion, we considered $K_1=K_2=51$ and we set $L=12$ to cut the tales of the unknown joint density of the underlying assets. 

\begin{table}[h]
\begin {scriptsize}
\begin{center}
\caption{Comparison between COS and finite volume IMEX Runge-Kutta numerical solutions for Test 1 at $t=T$. At each cell of the table the COS price is the first one, while the finite volume price is the second one.}\label{tb-basket-convective_anhelada}%
\begin{tabular}{c|c|c|c|c}
\backslashbox{$s_1$}{$s_2$} & $20.03125$ &$40.03125$ & $60.03125$ & $80.03125$  \\
\hline
 \multirow{2}{*}{$20.03125$} & $1.1495452217\times 10^{-11}$ & $3.5571109320$ & $13.5563429224$  & $23.5563429224$  \\
& $6.1707381613\times 10^{-12}$ &$3.5571143743$ & $13.5563429224$  & $23.5563429224$  \\
\hline
\multirow{2}{*}{$40.03125$} &$3.5571109320$& $13.5563429224$ & $23.5563429224$  & $33.5563429224$ \\
& $3.5571143743$ & $13.5563429224$ & $23.5563429224$  & $33.5563429224$ \\
\hline
\multirow{2}{*}{$60.03125$} &$13.5563429224$ & $23.5563429224$ & $33.5563429224$  & $43.5563429224$  \\
& $13.5563429224$ &$23.5563429224$ & $33.5563429224$  & $43.5563429224$  \\
\hline
\multirow{2}{*}{$80.03125$} &$23.5563429224$& $33.5563429224$ & $43.5563429224$ & $53.5563429224 $ \\
& $23.5563429224$ &$33.5563429224$ & $43.5563429224$  & $53.5563429224$  \\
\end{tabular}
\end{center}
\end{scriptsize}
\end{table}

\begin{table}[h]
\begin {scriptsize}
\begin{center}
\caption{Comparison between COS and finite volume IMEX Runge-Kutta numerical solutions for Test 2 at $t=T$. At each cell of the table the COS price is the first one, while the finite volume price is the second one.}\label{tb-basket-diffusive_anhelada}%
\begin{tabular}{c|c|c|c|c}
\backslashbox{$s_1$}{$s_2$} & $20.296875$ &$39.046875$ & $57.796875$ & $76.546875$  \\
\hline
 \multirow{2}{*}{$20.296875$} & $9.9328285143\times 10^{-2}$ & $2.7953041634$ & $10.1409569225$  & $19.2023244718$  \\
& $9.9312989509\times 10^{-2}$ &$2.7953229793$ & $10.1409654135$  & $19.2023248970$  \\
\hline
\multirow{2}{*}{$39.046875$} &$2.7953041634$& $10.0983319903$ & $19.1914361881$  & $28.5400105675$ \\
& $2.7953229793$ & $10.0983391832$ & $19.1914362212$  & $28.5400103410$ \\
\hline
\multirow{2}{*}{$57.796875$} &$10.1409569225$ & $19.1914361881$ & $28.5396217167$  & $37.9127244740$  \\
& $10.1409654135$ &$19.1914362212$ & $28.5396214779$  & $37.9127244317$  \\
\hline
\multirow{2}{*}{$76.546875$} &$19.2023244718$& $28.5400105675$ & $37.9127244740$ & $47.2875872432 $ \\
& $19.2023248970$ &$28.5400103410$ & $37.9127244317$  & $47.2875872518$  \\
\end{tabular}
\end{center}
\end{scriptsize}
\end{table}

\clearpage
\subsection{Heston model\label{sec:NumExpHeston}}

In this experiment we price a vanilla call option under the Heston stochastic volatility model. The PDE model \eqref{eq:HestonPDEForward} can be written in the conservative form \eqref{SistCons} as follows:
\begin{equation} \label{SistConsHeston}
    \dfrac{\partial u}{\partial t}  + \dfrac{\partial f_1}{\partial s} (u)+\dfrac{\partial f_2}{\partial v}  (u)=\dfrac{\partial g_1}{\partial s} (u_{s},u_{v})+\dfrac{\partial g_2}{\partial v} (u_{s},u_{v})+h(u),
\end{equation}
where the functions $f_1, f_2, g_1, g_2$ and $h$ are given by:
\begin{align*}
f_1(u)&=(v-r+q)s u(s,v,t),\\
f_2(u)&=\left(\rho \sigma v - \kappa(\theta-v) + \frac12 \sigma^2 \right)u(s,v,t),\\
g_1(u_s,u_v)&=\dfrac{1}{2}s^2v u_s(s,v,t)+\rho\sigma s v u_v(s,v,t),\\
g_2(u_s,u_v)&=\dfrac{1}{2}\sigma^2v u_v(s,v,t),\\
h(u)&=(v-2r+q+\kappa+\rho\sigma) u(s,v,t).
\end{align*}
The initial condition $u_0$, being equation \eqref{eq:hestonPayoff} for the European call option, was averaged on each volume. The initial volume averages are set as
$${\bar{u}^0}_{i,j} = \dfrac{1}{\lvert V_{ij} \rvert} \int_{V_{ij}} u_0(s)\dx{s}\dx{v}.$$

In this part we perform two vanilla call option pricing tests with the market parameters of Table \ref{tb-Heston}. The first set of parameters, labelled as Test 3, is taken from \cite{Foulon10}. The second group of parameters, denoted as Test 4, is a variation of Test 3, swapping the interest rate $r$ and the volatility of the volatility $\sigma$. In both tests, in order to minimize the numerical errors coming from the artificial boundary conditions, the computational space domain was considered as $(s,v)\in[0,800]\times[0,4]$.

\begin{table}[!htb]
\begin{center}
\caption{Market data for vanilla call options under Heston model.}
\begin{tabular}{c|c|c|c|c|c|c|c|c}
  &$\sigma$ & $r$ & $q$  & $\kappa$ & $\theta$  & $\rho$ & $T$& $K$\\
\hline
Test 3& $0.3$ & $0.025$  & \multirow{2}{*}{$0.0$} & \multirow{2}{*}{$1.5$} & \multirow{2}{*}{$0.04$} & \multirow{2}{*}{$-0.9$}& \multirow{2}{*}{$0.25$} & \multirow{2}{*}{$100$}\\

Test 4 &$0.025$ & $0.3$  &  &  &  &  &  & \\
\end{tabular} 
\label{tb-Heston}
\end{center}
\end{table}

\begin{figure}[!htb]
\centering
\subfigure[Test 1] {\includegraphics[height=4.4cm]{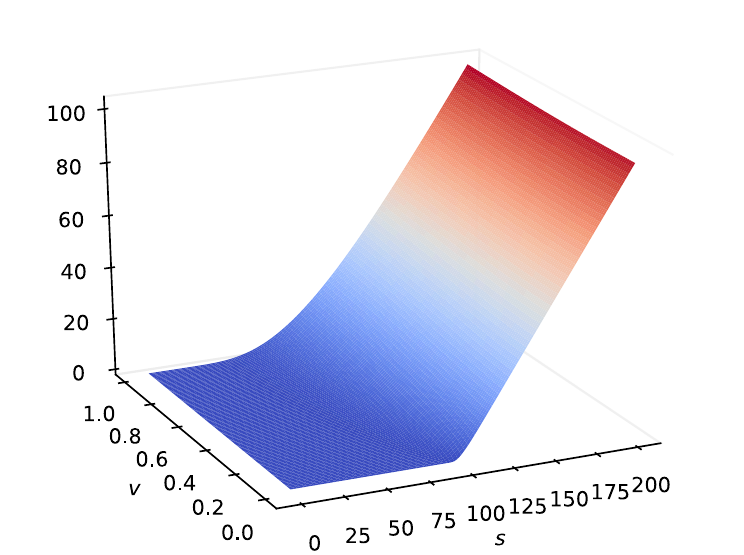}}
\subfigure[Test 2] {\includegraphics[height=4.4cm]{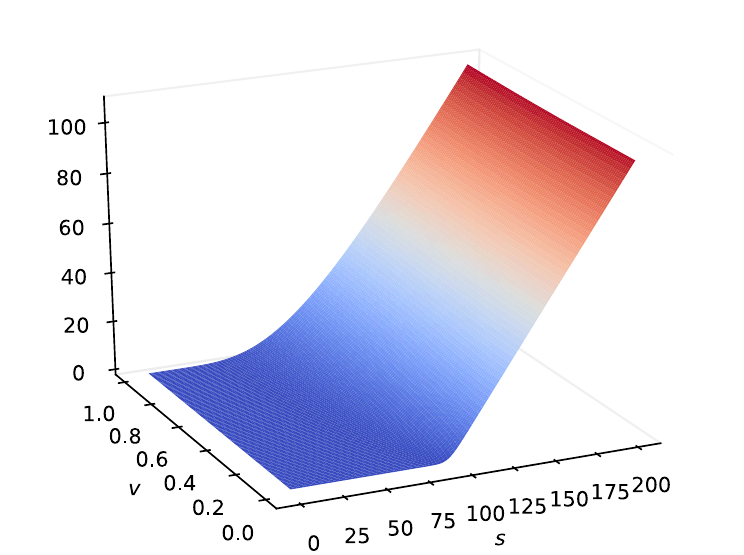}}
\caption{Vanilla call option prices for Test 3 (left) and Test 4 (right) at $t=T$.}
\label{fig:heston-prices}
\end{figure}

\begin{figure}[!htb]
\centering
\subfigure[Test 1] {\includegraphics[height=4.4cm]{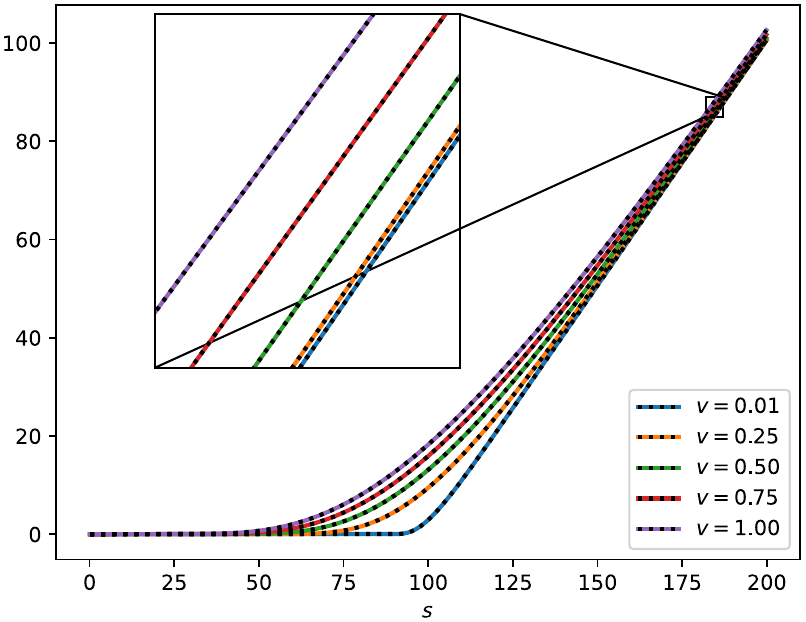}}
\subfigure[Test 2] {\includegraphics[height=4.4cm]{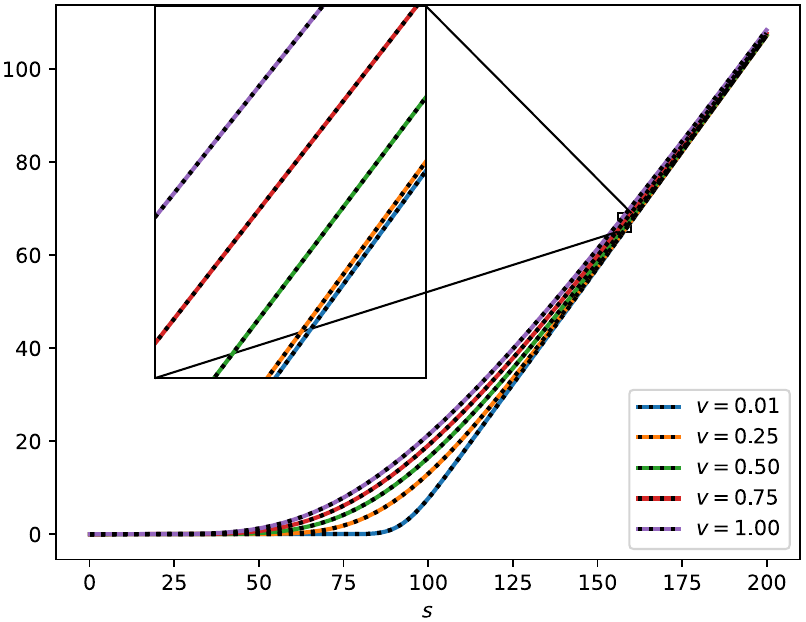}}
\caption{Cuts of prices surfaces of Figure \ref{fig:heston-prices}. Numerical solution is shown in continuous line and COS solution with squares.}
\label{fig:heston-prices-cuts}
\end{figure}

\begin{figure}[!htb]
\centering
\subfigure[Test 1] {\includegraphics[height=4.4cm]{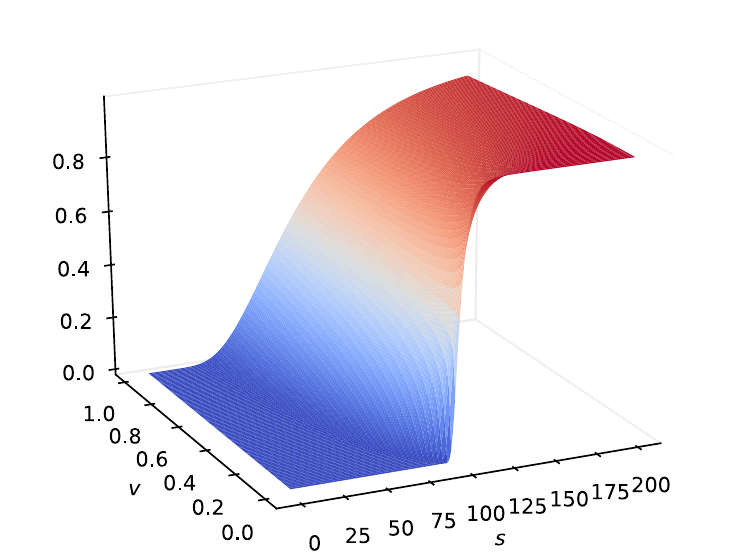}}
\subfigure[Test 2] {\includegraphics[height=4.4cm]{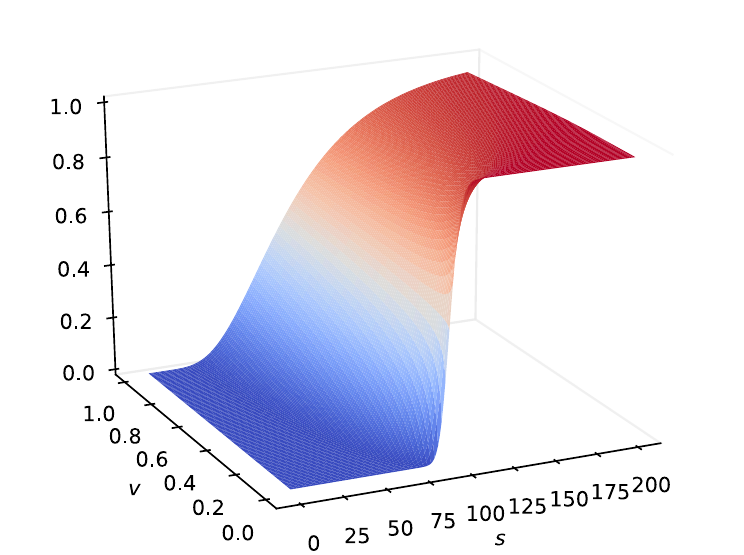}}
\caption{Vanilla call option deltas for Test 3 (left) and Test 4 (right) at $t=T$.}
\label{fig:heston-deltas}
\end{figure}

\begin{figure}[!htb]
\centering
\subfigure[Test 1] {\includegraphics[height=4.4cm]{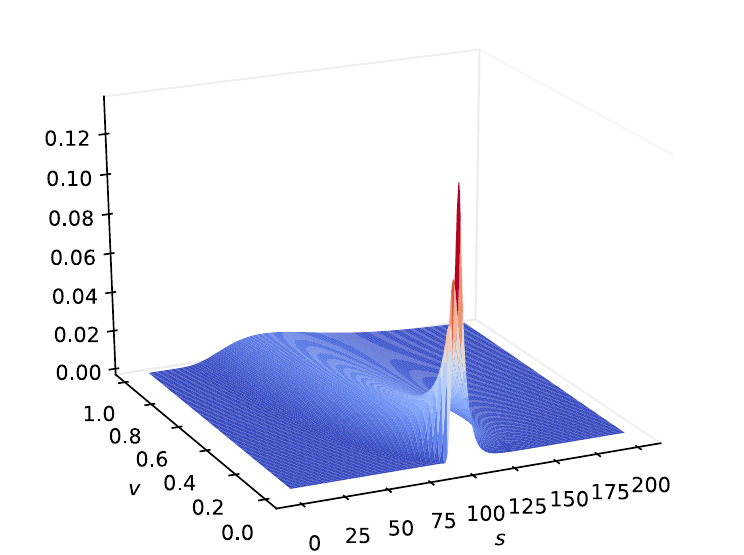}}
\subfigure[Test 2] {\includegraphics[height=4.4cm]{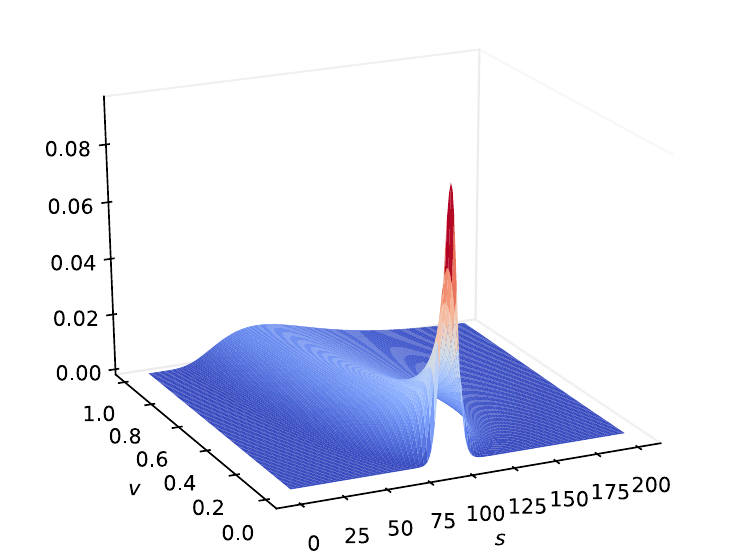}}
\caption{Vanilla call option gammas for Test 3 (left) and Test 4 (right) at $t=T$.}
\label{fig:heston-gammas}
\end{figure}

Figure \ref{fig:heston-prices} displays the finite volume IMEX Runge-Kutta approximations of the exact option price functions $u$ at $t=T$ corresponding to the two cases of Table \ref{tb-Heston}. Cuts of these price surfaces with planes parallel to $v=0$ are shown in Figure \ref{fig:heston-prices-cuts} in combination with COS solutions. Additionally, Figures \ref{fig:heston-deltas} and \ref{fig:heston-gammas} present the numerical Greeks of the numerical approximation, deltas and gammas, respectively. Once again, the here proposed numerical scheme provides accurate and non-oscillatory approximations of the Greeks, even at regions of discontinuities and non-smoothness in the initial condition.

At this stage, using the described COS method, we compute reference call option prices. For a mesh of size $3200\times 3200$, the errors of the approximation obtained with the finite volume IMEX Runge-Kutta method are presented in Table \ref{tb-heston_errors}. Figure \ref{tb-heston_errors} displays contour plots for the surface defined by the absolute errors in the call option prices at $t=T$. These plots show that the here developed finite volume schemes provide high-quality approximations even at the region of non-smoothness in the payoff.

\begin{table}[!htb]
\begin{footnotesize}
\begin{center}
\caption{Finite volume IMEX Runge-Kutta numerical errors against reference call option prices computed with COS method at $t=T$.}
\begin{tabular}{|c|c|c|}
\hline
\multicolumn{3}{|c|}{Test 3}\\
\hline
  $L_\infty  \mbox{ error }$ & $L_\infty \mbox{ relative error }$ & $\mbox{Mean absolute error}$ \\
\hline
$2.0377\times 10^{-4}$ & $3.2852\times 10^{-7}$ & $3.1711\times 10^{-5}$  \\
\hline

\hline
\multicolumn{3}{|c|}{Test 4}\\
\hline
  $L_\infty  \mbox{ error }$ & $L_\infty \mbox{ relative error }$ & $\mbox{Mean absolute error}$ \\
\hline
$1.5196\times 10^{-4}$ & $1.9969\times 10^{-5}$ & $2.6785\times 10^{-7}$  \\
\hline
\end{tabular} 
\label{tb-heston_errors} 
\end{center}
\end{footnotesize}
\end{table}

\begin{figure}[!htb]
\centering
\subfigure[Test 1] {\includegraphics[height=4.4cm]{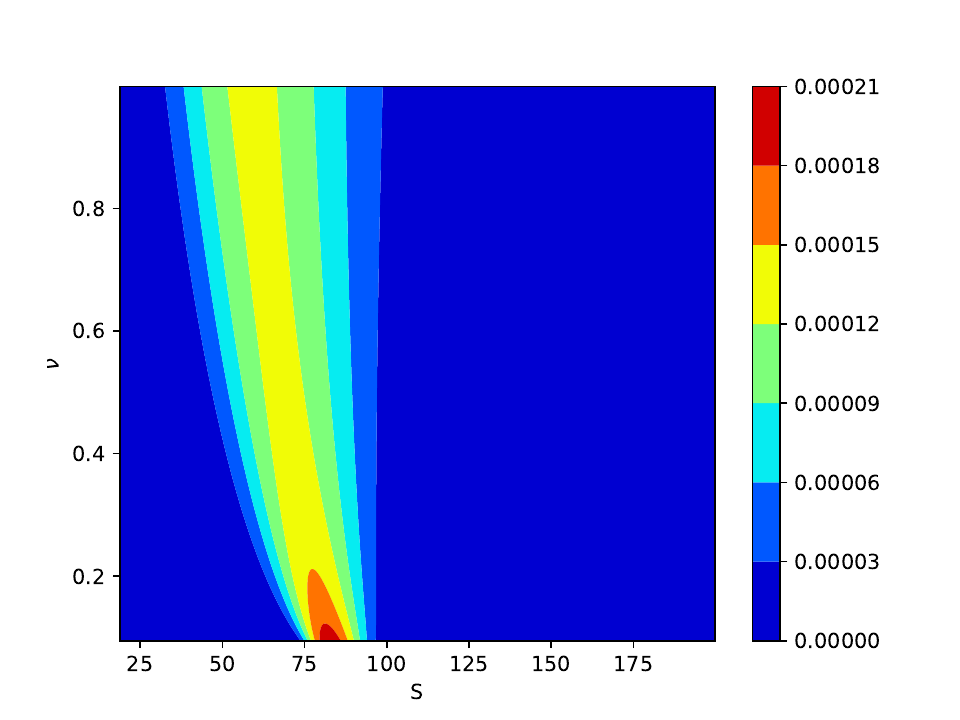}}
\subfigure[Test 2] {\includegraphics[height=4.4cm]{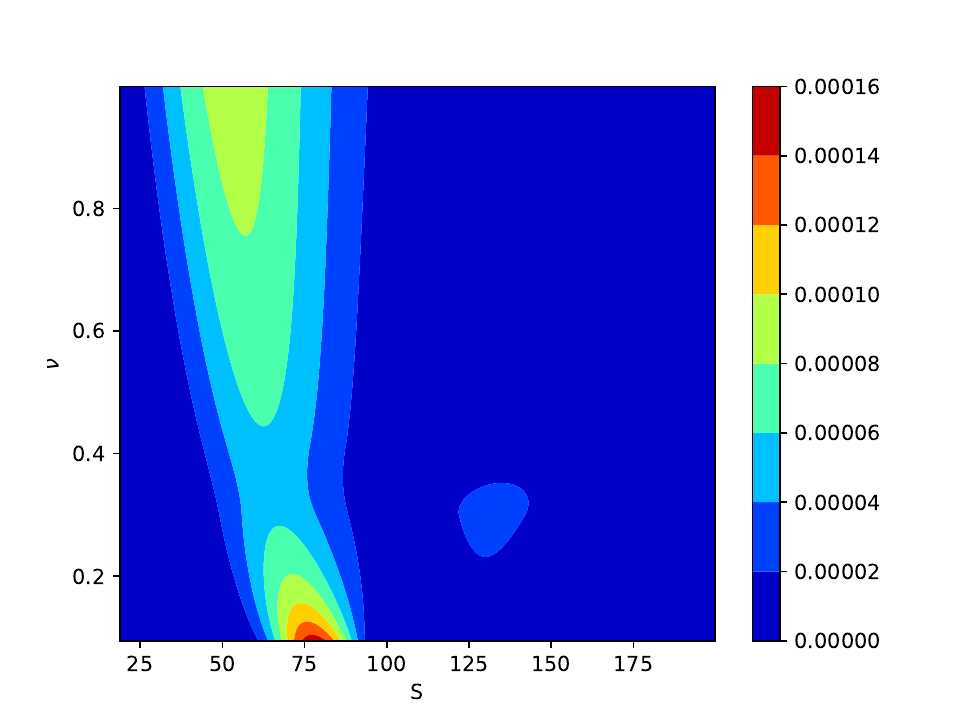}}
\caption{Contour plots for the surface of absolute errors of call option prices for Test 3 (left) and Test 4 (right) at $t=T$.}
\label{fig:error_heston_65}
\end{figure}

In Tables \ref{tb-Heston65} and \ref{tb-Heston66}, the $L_1$ errors and the $L_1$ orders of the explicit and IMEX methods for the Test 3 and Test 4 are displayed, respectively. Again, both numerical schemes achieve second-order accuracy in the $L_1$ norm. Even if we impose time-dependent boundary conditions, we do not observe the troublesome order reduction effect described in other papers. As said previously, IMEX is able to converge using much larger time steps than the explicit method, thus it consumes much less computing time.

\begin{table}[!h]
\begin{center}
{\scriptsize
\begin{tabular}{|c||c|c|c|c|}
\hline
& \multicolumn{4}{c|}{IMEX} \\
\hline
$N_1\times N_2$  & $L_1$ \text{ error} &  \text{ Order}& $\Delta t$ & Time (s) \\
\hline
$25\times 25$ &  $3.1356\times 10^{2}$ & $--$ & $3.85\times 10^{-3}$ & $1.6 \times 10^{-2}$ \\
$50\times 50$ &  $1.0516\times 10^{2}$ & $1.58$&  $1.92\times 10^{-3}$ & $9.7\times 10^{-2}$\\
$100\times 100$ &  $2.5715\times 10^{1}$ & $2.03$&  $9.61\times 10^{-4}$ & $7.4 \times 10^{-1}$ \\
$200\times 200$ &  $6.2028\times 10^{0}$ & $2.05$ &  $4.80\times 10^{-4}$ & $7.4 \times 10^{0}$ \\
$400\times 400$ &  $1.4586\times 10^{0}$ & $2.09$ &  $2.40\times 10^{-4}$ & $6.8 \times 10^{1}$ \\
$800\times 800$ &  $2.9016\times 10^{-1}$ & $2.33$ &  $1.20\times 10^{-4}$ & $6.4 \times 10^{2}$ \\
\hline
\hline
& \multicolumn{4}{c|}{Explicit} \\
\hline
$N_1\times N_2$  & $L_1$ \text{ error} &  \text{ Order}& $\Delta t$ & Time (s) \\
\hline
$25\times 25$ & $3.1312\times 10^{2}$ & $--$ & $1.92\times 10^{-4}$  & $1.9 \times 10^{-1}$\\
$50\times 50$ & $1.0489\times 10^{2}$ &$1.58$ &  $4.81\times 10^{-5}$ & $ 2.7\times 10^{0}$ \\
$100\times 100$ & $2.5649\times 10^{1}$ &$2.03$ &  $1.20\times 10^{-5}$ & $4.4 \times 10^{1}$ \\
$200\times 200$ & $6.1864\times 10^{0}$ &$2.05$ &  $3.00\times 10^{-6}$ & $7.3 \times 10^{2}$ \\
$400\times 400$ & $1.4543\times 10^{0}$ &$2.09$ &  $7.52\times 10^{-7}$ & $1.2 \times 10^{4}$ \\
$800\times 800$ & $2.8924\times 10^{-1}$ &$2.33$ &  $1.87\times 10^{-7}$ & $2.0 \times 10^{5}$ \\
\hline
\end{tabular} 
}
\end{center}
\caption{$L_1$ errors and $L_1$ orders of convergence of the IMEX and explicit finite volume methods for Test 3.} 
\label{tb-Heston65} 
\end{table}

\begin{table}[!h]
\begin{center}
{\scriptsize
\begin{tabular}{|c||c|c|c|c|}
\hline
& \multicolumn{4}{c|}{IMEX} \\
\hline
$N_1\times N_2$  & $L_1$ \text{ error} &  \text{ Order}& $\Delta t$ & Time (s)  \\
\hline
$25\times 25$ &  $8.9656\times 10^{1}$ & $--$ & $3.87\times 10^{-3}$ & $1.5\times 10^{-2}$ \\
$50\times 50$ &  $2.4203\times 10^{1}$ & $1.89$& $1.94\times 10^{-3}$ & $9.0\times 10^{-2}$\\
$100\times 100$ &  $9.3022\times 10^{0}$ & $1.38$&  $9.69\times 10^{-4}$ & $7.0\times 10^{-1}$ \\
$200\times 200$ &  $1.2578\times 10^{0}$ & $2.89$ & $4.84\times 10^{-4}$ & $6.9\times 10^{0}$ \\
$400\times 400$ &  $2.9883\times 10^{-1}$ & $2.07$ & $2.42\times 10^{-4}$  & $6.8\times 10^{1}$ \\
$800\times 800$ &  $5.9846\times 10^{-2}$ & $2.32$ & $1.21\times 10^{-4}$ & $6.3\times 10^{2}$ \\
\hline
\hline
& \multicolumn{4}{c|}{Explicit} \\
\hline
$N_1\times N_2$  & $L_1$ \text{ error} &  \text{ Order}& $\Delta t$ & Time (s)  \\
\hline
$25\times 25$ & $9.0576\times 10^{1}$ & $--$  & $1.99\times 10^{-4}$ & $1.8\times 10^{-1}$\\
$50\times 50$ & $2.4440\times 10^{1}$ &$1.89$  & $4.99\times 10^{-5}$ & $2.7\times 10^{0}$\\
$100\times 100$ & $9.3648\times 10^{0}$ &$1.38$ & $1.25\times 10^{-5}$ & $4.4\times 10^{2}$ \\
$200\times 200$ & $1.2671\times 10^{0}$ &$2.89$ & $3.12\times 10^{-6}$ & $7.1\times 10^{2}$ \\
$400\times 400$ & $3.0089\times 10^{-1}$ &$2.07$ & $7.79\times 10^{-7}$ & $1.2\times 10^{4}$ \\
$800\times 800$ & $6.0411\times 10^{-2}$ &$2.32$  &  $1.95\times 10^{-7}$ & $1.9\times 10^{5}$ \\
\hline
\end{tabular} 
}
\end{center}
\caption{$L_1$ errors and $L_1$ orders of convergence of the IMEX and explicit finite volume methods for Test 4.}
\label{tb-Heston66} 
\end{table}

Finally, Tables \ref{tb-Heston65-anhelada} and \ref{tb-Heston65-anheladabis} show some vanilla call option prices at $t=T$ under Heston stochastic volatility model, for Test 3 and Test 4, respectively. These prices were obtained with COS method and with the finite volume IMEX Runge-Kutta solver executed over a grid with size $3200\times 3200$. 

\begin{table}[!h]
\begin{scriptsize}
\begin{center}
\caption{Comparison between COS and finite volume IMEX Runge-Kutta numerical solutions for Test 3 at $t=T$. At each cell of the table the COS price is the first one, while the finite volume price is the second one.}
\begin{tabular}{c|c|c|c|c}
 \backslashbox{$v$}{$s$}  & $75.125$ & $100.125$ & $125.125$ & $150.125$ \\
\hline
$0.200625$ & $0.4317502658$ & $8.5901684247$ & $27.6694906563$ &  $51.1935422402$\\
 &  $0.4316035999$ & $8.5901562104$ & $27.6695002425$ &  $51.1935390259$\\
 \hline
$0.400625$ & $1.8664224003$ & $11.8552660514$ & $30.0081668834$ & $52.3493028388$ \\
 &  $1.8662948639$ & $11.8552481800$ & $30.0081846067$ &  $52.3493135590$\\
 \hline
$0.600625$ & $3.3658574222$ & $14.3630707461$ & $32.1382630288$ & $53.7378969543$  \\
 &  $3.3657449349$ & $14.3630493039$ & $32.1382791463$ &  $53.7379138145$\\
 \hline
$0.800625$ & $4.7800189282$ & $16.4716211703$ & $34.0691648265$ & $55.1769879733$ \\
 &  $4.7799153837$ & $16.4715967918$ & $34.0691774343$ &  $55.1770065036$\\
\end{tabular} 
\label{tb-Heston65-anhelada}
\end{center}
\end{scriptsize}
\end{table}

\begin{table}[!h]
\begin{scriptsize}
\begin{center}
\caption{Comparison between COS and finite volume IMEX Runge-Kutta numerical solutions for Test 4 at $t=T$. At each cell of the table the COS price is the first one, while the finite volume price is the second one.}\label{tb-Heston65-anheladabis}%
\begin{tabular}{c|c|c|c|c}
\backslashbox{$v$}{$s$} & $75.125$ &$100.125$ & $125.125$ & $150.125$  \\
\hline
 \multirow{2}{*}{$0.200625$} & $1.3840492642$ & $12.2239735395$ & $33.1373172158$  & $57.4478630579$  \\
& $1.3839721771$ & $12.2239654379$  & $33.1372986255$   & $57.4478503229$   \\
\hline
\multirow{2}{*}{$0.400625$} & $3.2984130167$ & $15.2478827153$  & $34.8602149676$   & $58.0809835379$  \\
& $3.2983708067$  & $15.2478733955$  & $34.8602005816$   & $58.0809679113$  \\
\hline
\multirow{2}{*}{$0.600625$} & $5.0070610101$  & $17.6182306257$  & $36.6338455366$   & $59.0749794899$   \\
& $5.0070098251$  & $17.6182226102$  & $36.6338456927$   & $59.0749764400$   \\
\hline
\multirow{2}{*}{$0.800625$} & $6.5349321045$  & $19.6278242193$  & $38.3205477349$  & $60.2165220825$  \\
& $6.5348770744$  &$19.6278146746$  & $38.3205529668$   & $60.2165267905$  \\
\end{tabular}
\end{center}
\end{scriptsize}
\end{table}

\clearpage
\section{Conclusions} \label{conclusions}

We have shown that finite volume IMEX Runge-Kutta numerical schemes are suitable for solving convection-diffusion PDE option pricing problems, even in the presence of mixed derivatives in the diffusive operator. This opens the door to the application of these schemes to numerous models in finance, even those giving rise to non-linear PDEs.

The obtained numerical schemes are highly efficient. On the one hand, large time steps can be used, avoiding the need to use small time steps enforced by the diffusion stability condition that appears when explicit schemes are considered. On the other hand, the schemes are second order accurate. This fact is of paramount importance in order to obtain accurate approximations of the Greeks without oscillations. Besides, we have shown that second order accuracy is preserved even when non smooth initial conditions (payoffs) are considered, which is the usual situation in finance. Thus, additional special techniques, like the Rannacher time stepping, do not need to be taken into account. Moreover, the here developed option price calculators can be extended to build numerical solvers with higher order convergence.

Last, but not least, in this work we provide an alternative way to compute very accurate approximations of the prices of arithmetic basket call options by means of a highly accurate and efficient multidimensional COS Fourier method. These Fourier semi-analytical solutions can be also very valuable benchmarks for a broader audience that work in the development of high order schemes in the general parabolic setting, even outside the financial world.

\bibliography{mybibfile}

\end{document}